\documentclass[12pt]{amsart} 

\begin{document}

%%%%%%%%%%%%%%%%%%%%%%%%%%%%%%%%%%%%%%%%%%%%%%%%%%%%%%
%%%%%%%%%%%%%%%%%%%%%%%%%%%%%%%%%%%%%%%%%%%%%%%%%%%%%%
\renewcommand{\theequation}{\thesection.\arabic{equation}}
\def\theenumi{\roman{enumi}}
\def\labelenumi{(\theenumi)}
%%%%%%%%%%%%%%%%%%%%%%%%%%%%%%%%%%%%%%%%%%%%%%%%%%%%%%%
%loading  gothic fonts
\font\teneufm=eufm10
\font\seveneufm=eufm7
\font\fiveeufm=eufm5
\newfam\eufmfam
\textfont\eufmfam=\teneufm
\scriptfont\eufmfam=\seveneufm
\scriptscriptfont\eufmfam=\fiveeufm
\def\frak#1{{\fam\eufmfam\relax#1}}
\let\goth\frak
%%%%%%%%%%%%%%%%%%%%%%%%%%%%%%%%%%%%%%%%%%%%%%%%%%%%%
\newcommand{\g}{{\goth{g}}}
\newcommand{\slth}{\widehat{\mbox{\twelveeufm sl}}_2} %??
\newcommand{\slN}{\mbox{\twelveeufm sl}_N} %??
\newcommand{\slthten}{\widehat{\mbox{\teneufm sl}}_2} %??
%%%%%%%%%%%%%%%%%%%%%%%%%%%%%%%%%%%%%%%%%%%%%%%%%%%%%
\font\twelveeufm=eufm10 scaled\magstep1
\font\fourteeneufm=eufm10 scaled\magstep2    %for title,subsection
\font\seventeeneufm=eufm10 scaled\magstep3   %for section
\font\twentyoneeufm=eufm10 scaled\magstep4   
\newcommand{\slNBig}{\mbox{\seventeeneufm sl}_N} %??
%\newcommand{\slNBig}{\goth{sl}_N} %???
%%%%%%%%%%%%%%%%%%%%%%%%%%%%%%%%%%%%%%%%%%%%%%%%%%%%%
\newcommand{\Z}{{\mathbb Z}} %??
\newcommand{\C}{{\mathbb C}} %??
\newcommand{\R}{{\mathbb R}} %??

\newcommand{\F}{{\mathcal F}}
\newcommand{\Ft}{\widetilde{\mathcal F}}
\renewcommand{\H}{{\mathcal H}}
\newcommand{\Cc}{{\mathcal C}}
\newcommand{\ve}{\varepsilon}
\newcommand{\la}{\lambda}
\newcommand{\Wb}{\overline{W}}
\newcommand{\dz}{\underline{dz}}
\newcommand{\dw}{\underline{dw}}
\newcommand{\dbr}[1]{{\langle\!\langle #1 \rangle\!\rangle}} %double bracket
\newcommand{\G}[3]{{{G}_{x^{#1}}(#2,#3)}}
\newcommand{\dsp}{\displaystyle}
\newcommand{\scr}{\scriptstyle}
\newcommand{\qi}[1]{{ [  #1 ]_x}}

\newcommand{\nn}{\nonumber}
\newcommand{\bea}{\begin{eqnarray}}
\newcommand{\ena}{\end{eqnarray}}
\newcommand{\be}{\begin{eqnarray*}}
\newcommand{\en}{\end{eqnarray*}}
\newcommand{\lb}[1]{\label{#1}}

\newcommand{\res}{{\mathop{\rm res}}}
\newcommand{\id}{{\rm id}}
\newcommand{\tr}{{\rm tr}}
\newcommand{\bra}[1]{\langle #1 |}        %bra
\newcommand{\ket}[1]{{| #1 \rangle}}      %ket
\newcommand{\rs}{{r^*}}
\newcommand{\BW}[5]{\Bigl({#1\atop#3}\ {#2\atop#4}\Bigl|#5\Bigr)}
\newcommand{\sfrac}[2]{{\textstyle\frac{#1}{#2}}}
%%%%%%%%%%%%%%%%%%%%%%%%%%%%%%%%%%%%%%%%%%%%%%%%%%%%%%%%%%
\newcommand{\Remark}{\medskip \noindent {\it Remark.}\quad}
\newcommand{\example}{\medskip \noindent {\it Example.}\quad}

\newtheorem{thm}{Theorem}[section]
\newtheorem{conj}[thm]{Conjecture}
\newtheorem{prop}[thm]{Proposition}
\newtheorem{lem}[thm]{Lemma}
\newtheorem{dfn}[thm]{Definition}
\newcommand{\ignore}[1]{}
%%%%%%%%%%%%%%%%%%%%%%%%%%%%%%%%%%%%%%%%%%%%%%%%%%%%%%%%%%
%%%%%%%%%%%%%%%%%%%%%%%%%%%%%%%%%%%%%%%%%%%%%%%%%%%%%%%%%%

\title{Free Field Construction for the ABF Models in Regime II}

\author[M.~Jimbo]{Michio~Jimbo}
\address{Division of Mathematics, Graduate School of Science\\
Kyoto University, Kyoto 606-8502, Japan}
\email{jimbo@kusm.kyoto-u.ac.jp}
\author[H.~Konno]{Hitoshi~Konno}
\address{Department of Mathematics, Faculty of Integrated Arts and Sciences\\
Hiroshima University, Higashi-Hiroshima 739-8521, Japan}
\email{konno@mis.hiroshima-u.ac.jp}	
\author[S.~Odake]{Satoru~Odake}
\address{Department of Physics, Faculty of Science\\
Shinshu University, Matsumoto 390-8621, Japan}
\email{odake@azusa.shinshu-u.ac.jp}
\author[Y.~Pugai]{Yaroslav~Pugai}
\address{Research Institute for Mathematical Sciences\\
Kyoto University, Kyoto 606-8502, Japan}
\email{slava@kurims.kyoto-u.ac.jp}
\author[J.~Shiraishi]{Jun'ichi~Shiraishi}
\address{Graduate School of Mathematical Sciences\\
The University of Tokyo, Tokyo 153-8914, Japan}
\email{shiraish@ms.u-tokyo.ac.jp}

\date{\today}
\dedicatory{Dedicated to Professor Rodney Baxter on the occasion
of his sixtieth birthday}

\setcounter{footnote}{0}\renewcommand{\thefootnote}{\arabic{footnote}}

\begin{abstract}
The Wakimoto construction for the quantum affine algebra $U_q(\slthten)$ 
admits a reduction to the $q$-deformed parafermion algebras. 
We interpret the latter theory as a free field realization of 
the Andrews-Baxter-Forrester models in regime II. 
We give multi-particle form factors of some local operators  
on the lattice and compute their scaling limit, where 
the models are described by a massive field theory 
with $\Z_k$ symmetric minimal scattering matrices. 
\end{abstract}

\maketitle

\vskip20mm
\noindent
math.QA/0001071

\vskip3mm\noindent
HU-IAS/K-8\\
DPSU-99-8\\
RIMS-1266

\newpage

\setcounter{footnote}{0}
\renewcommand{\thefootnote}{\arabic{footnote})}
\renewcommand{\arraystretch}{1.2}

\setcounter{section}{0}
\setcounter{equation}{0}
%%%%%%%%%%%%%%%%%%%%%%%%%%%%%%%%%%%%%%%%%%%%%%%%%%%%%%%%%%%%%%%
%                                                             %
%  1. Introduction                                            %
%                                                             %
%%%%%%%%%%%%%%%%%%%%%%%%%%%%%%%%%%%%%%%%%%%%%%%%%%%%%%%%%%%%%%%
\section{Introduction} \lb{sec:1}

In conformal field theory (CFT), 
free field construction (or the \linebreak
`Coulomb~gas' representation) \cite{DoFa84}
is the most effective calculational tool for physical quantities.  
The same can be said about
the vertex operator approach to solvable lattice models (see e.g. \cite{JM}). 
The latter can be viewed as a $q$-deformation of 
some conformal field models,   
the most typical example being 
the Andrews-Baxter-Forrester (ABF) models in regime III \cite{LP96} 
which corresponds to the minimal unitary series. 
Despite its technical importance, however, 
no uniform recipe is known at present for finding
a free field representation of a given CFT.
Even when it is known, an equally non-trivial task is to 
do the same for the corresponding off-critical solvable lattice models. 
In this article we address this issue 
in the case of the ABF models in regime II \cite{ABF}.  

By the level-rank duality for the Boltzmann weights \cite{JMO1}, 
the ABF models in regime II at level $k$ are equivalent to 
the $A^{(1)}_{k-1}$ face models in regime III at level $2$.  
For the latter, a free field realization is already available \cite{AJMP}. 
However this picture is rather complicated, since 
one has to deal with $k-1$ kinds of oscillators. 
The Felder complex in higher rank is also quite cumbersome \cite{FJMOP}. 
On the other hand, the ABF models in regime II 
are described in the conformal limit 
by the $\Z_k$ parafermionic CFT \cite{ZaFa85},  
and in the massive scaling limit 
by its perturbation by the first energy operator 
\cite{Za88,Ts88,BaRe90,Fat91,KlMe90}. 
Therefore one naturally expects that 
an alternative construction on the lattice 
is to invoke the $q$-deformation of the 
parafermion theory \cite{Ma94,Konno},     
where it is sufficient to deal with only two kinds of oscillators 
and the resolution of Fock modules has a simpler structure.  
In this paper we work out this point in detail. 
Our purpose is to show that the 
known results in representation theory fits nicely 
the interpretation as a bosonization of lattice models.  

The text is organized as follows. 
In section 2 we recall the ABF models and set up the notation.  
In section 3 we review the free field realization of the 
deformed parafermion theory. 
We then present various commutation relations of operators. 
Comparing them with those in the lattice theory,  	
we identify the deformed counterpart of the 
primary fields and the parafermionic currents  
with the vertex operators (VO's) of type I and type II, respectively. 
A simplifying feature is that 
these deformed parafermionic currents mutually commute by scalars. 
To our knowledge, this is the only known example 
in the vertex operator approach\footnote{as well as the 
lattice counterpart of the non-unitary minimal models $M_{2,2n+1}$ 
which are expected to have the same property.}
where the type II VO's do not involve integrals of screening operators. 
We note also the existence of a deformation of the $W_k$ currents, 
which we believe to be equivalent to the one in \cite{AKOS96,FeFr95}
with specialization $r=k+2$.
This conforms with the known equivalence between 
the $\Z_k$ parafermionic CFT 
and the first member in the $W_k$ minimal series.  
In section 4 we calculate form factors of height variables on one 
and neighboring two lattice sites. 
In section 5 we take their continuous limit. 
Section 6 is devoted to the summary and discussions 
about open problems. 
Some technical formulas are gathered in the Appendices.

\setcounter{section}{1}
\setcounter{equation}{0}
%%%%%%%%%%%%%%%%%%%%%%%%%%%%%%%%%%%%%%%%%%%%%%%%%%%%%%%%%%%%%%%
%                                                             %
%  2. ABF model in regime II                                  %
%                                                             %
%%%%%%%%%%%%%%%%%%%%%%%%%%%%%%%%%%%%%%%%%%%%%%%%%%%%%%%%%%%%%%%
\section{ABF model in regime II}\lb{sec:2}

%%%%%%%%%%%%%%%%%%%%%%%%%%%%
%  2.1                     %
%%%%%%%%%%%%%%%%%%%%%%%%%%%%
\subsection{Boltzmann weights}\lb{subsec:2.1}

In this section we recall the ABF model in regime II 
and set up the notation. 
Throughout this article, we fix a positive integer $k\ge 2$. 

The Boltzmann weights of the ABF model \cite{ABF} have the form 
\bea
&&W\BW{a}{b}{c}{d}{u}= \rho(u)\Wb\BW{a}{b}{c}{d}{u}.   
\lb{BW0}
\ena
Here $a,b,c,d$ denote positive integers which we refer to as 
height variables.  
Besides, the weights depend on two real parameters $u$ and $x$. 	
If 
$1\le a,b,c,d\le k+1$ and $|a-b|=|b-d|=|d-c|=|c-a|=1$,  
then we set
\bea
&&\Wb\BW{a}{a\pm 1}{a\pm 1}{a\pm 2}{u}=1, 
\lb{BW1} \\
&&\Wb\BW{a}{a\pm 1}{a\pm 1}{a}{u}=\frac{[a\pm u][1]}{[a][1-u]},  
\lb{BW2}\\
&&\Wb\BW{a}{a\pm 1}{a\mp1}{a}{u}=\frac{[a\mp 1 ][-u]}{[a][1-u]}. 
\lb{BW3}
\ena
Here $[u]$ stands for the function 
\be
&&[u]= x^{\frac{u^2}{k+2}-u} \Theta_{x^{2(k+2)}}(x^{2u}),
\en
and we use the standard symbols 
\be
&&\Theta_{p}(z)=(z;p)_\infty (pz^{-1};p)_\infty (p;p)_\infty, 
\\
&&
(z_1,\cdots,z_n;p_1,\cdots, p_m)_\infty
=\prod_{j=1}^n\prod_{l_1,\cdots,l_m\ge 0}
(1-p_1^{l_1}\cdots p_m^{l_m}z_j).
\en
For all other values of $a,b,c,d$, 
we set the weight \eqref{BW0} to $0$.
 
In this note we restrict to regime II defined by 
\bea
&&0<x<1,\qquad -\frac{k}{2}<u<0.
%\lb{reg2}
\nn
\ena
We choose the scalar factor $\rho(u)$ in \eqref{BW0} to ensure that 
the partition function per site equals to $1$.   
Explicitly 
\begin{eqnarray*}
&&\rho(u)=x^{\frac{2u}{k(k+2)}}\frac{\rho_+(u)}{\rho_+(-u)},\\
&&\rho_+(u)=\frac{(x^{2k+2+2u},x^{2k+2+2u};x^{2k},x^{2(k+2)})_\infty}
{(x^{2k+2u},x^{2k+4+2u};x^{2k},x^{2(k+2)})_\infty}.
\end{eqnarray*}
We have 
\begin{eqnarray*}
&& \rho(u)\rho(-u)=1,\qquad
\rho(u)\rho(-k-u)=\frac{[1-u]^2}{[-u][2-u]}.
\end{eqnarray*}
In \cite{LP96}, a free field construction of the ABF model 
was found in regime III defined by $0<x<1, 0<u<1$.   
The Boltzmann weights in regime II differs from that case 
only in the choice of $\rho(u)$ made above. 

We represent the Boltzmann weights graphically as follows.
$$
\setlength{\unitlength}{1mm}
\begin{picture}(100,25)(10,8)
\put(25,10){\makebox(30,20)[l]{${\displaystyle W\BW{a}{b}{c}{d}{u}\;\;=}$}}
\put(60,14){\vector(1,0){7}}\put(67,14){\line(1,0){5}}
\put(60,26){\vector(1,0){7}}\put(67,26){\line(1,0){5}}
\put(60,26){\vector(0,-1){7}}\put(60,19){\line(0,-1){5}}
\put(72,26){\vector(0,-1){7}}\put(72,19){\line(0,-1){5}}
\put(59,27){\makebox(0,0)[br]{$a$}}
\put(59,13){\makebox(0,0)[tr]{$c$}}
\put(73,27){\makebox(0,0)[bl]{$b$}}
\put(73,13){\makebox(0,0)[tl]{$d$}}
\put(66,20){\makebox(0,0){$u$}}
\end{picture}
$$

%%%%%%%%%%%%%%%%%%%%%%%%%%%%
%  2.2                     %
%%%%%%%%%%%%%%%%%%%%%%%%%%%%
\subsection{Vertex operators in the na\"\i ve picture}\lb{subsec:2.2}

In this subsection we summarize the vertex operator approach 
in the na\"\i ve picture \cite{Fo94}.  

In regime II, there are $2k$ different ground states  
labeled by $m\in \Z/2k\Z$. 
Each of them is invarinat under the shift in the NE-SW direction, and hence 
is specified by a sequence of heights on a column. 
Fix a reference column and a site (say, site $1$) on it. 
Then the $m$-th ground state is given by the following 
sequence $\{\bar{l}^{(m)}_j\}_{j\in\Z}$ on that column
\be
&&\bar{l}^{(m)}_{j-m}=
\begin{cases}
j  & (1\le j\le k) \\
2k+2-j & (k+1\le j\le 2k)
\end{cases},\\
&&\bar{l}^{(m)}_{j+2k}=\bar{l}^{(m)}_{j} ~~(j\in\Z). 
\en
Here the sites on the reference column are numbered by $j$ from 
south to north in the increasing order. 

Let $\H_{m,a}$ 
denote the space of states of the half-infinite lattice, 
in the sector where the central height 
(i.e. the one on the reference site $1$) 
is fixed to $a$ and the boundary heights are in the ground state $m$. 
Namely we set formally 
\be
\H_{m,a}=\mbox{ span }\Bigl\{(l_j)_{j=1}^\infty
\mid 1\le l_j\le k+1,~~
l_1=a,~~l_j=\bar{l}^{(m)}_j~~(j\gg 1)\Bigr\}.
\en
Notice that $\H_{m,a}=0$ if $a\equiv m \bmod 2$. 
The corner transfer matrices (CTM's)  
associated with the north-west, south-west, south-east 
and north-east quadrants operate respectively as
\be
A^{(a)}_{NW}(u)&:& \H_{m,a} \longrightarrow \H_{m,a}, \\
A^{(a)}_{SW}(u)&:& \H_{m,a} \longrightarrow \H_{-m,a},\\
A^{(a)}_{SE}(u)&:& \H_{-m,a} \longrightarrow \H_{-m,a}, \\
A^{(a)}_{NE}(u)&:& \H_{-m,a} \longrightarrow \H_{m,a}. 
\en
In the limit of an infinitely large lattice, 
they have the simple form \cite{ABF}
\be
&&A^{(a)}_{NW}(u)=x^{-2uH^{(a)}_C},\\%\lb{CTM1}\\
&&A^{(a)}_{SW}(u)=M^{(a)}x^{2uH^{(a)}_C},\\%\lb{CTM2}\\
&&A^{(a)}_{SE}(u)=M^{(a)}x^{-2uH^{(a)}_C}{M^{(a)}}^{-1},\\%\lb{CTM3}\\
&&A^{(a)}_{NE}(u)=[a] x^{2(u+k)H^{(a)}_C}{M^{(a)}}^{-1},%\lb{CTM4}
\en
where the `corner Hamiltonian' $H^{(a)}_C:\H_{m,a}\rightarrow\H_{m,a}$ 
has a discrete and equidistant spectrum,  
and $M^{(a)}:\H_{m,a}\rightarrow\H_{-m,a}$ is an operator which 
does not play a role in the following. 
Both of them are independent of $u$. 

Likewise let 
\be
\Phi^{(a,b)}_N(u) &:& \H_{m,b} \longrightarrow \H_{m+1,a}, \\
\Phi^{(b,a)}_W(-u) &:& \H_{m+1,a} \longrightarrow \H_{m,b}, \\
\Phi^{(a,b)}_S(u)  &:& \H_{-m,b} \longrightarrow \H_{-m+1,a}, \\
\Phi^{(b,a)}_E(-u) &:& \H_{-m+1,a} \longrightarrow \H_{-m,b}
\en
be the half-infinite transfer matrices extending to 
north, west, south and east directions, respectively. 
(Note that the $(m+1)$-th ground state is obtained 
from the $m$-th one by shifting it one step to the left.)  
Normally we suppress the dependence on the boundary condition $m$ 
in the notation.  
Unlike the case of regime III, 
our Boltzmann weights do not have the crossing symmetry. 
Accordingly there is no simple relation between the operators 
$\Phi_N^{(a,b)}(u)$ and $\Phi_W^{(a,b)}(u)$.
\medskip

$$
\setlength{\unitlength}{1mm}
\begin{picture}(140,90)(-10,-7)
\put(50,50){
%%%%% begin subfig 
\begin{picture}(15,35)(10,10)
\put(10,38){\vector(0,-1){18.5}}
\put(10,19.5){\vector(0,-1){7}}\put(10,12.5){\line(0,-1){2.5}}
\put(17,38){\vector(0,-1){18.5}}
\put(17,19.5){\vector(0,-1){7}}\put(17,12.5){\line(0,-1){2.5}}
\put(10,10){\vector(1,0){4.5}}\put(14.5,10){\line(1,0){2.5}}
\put(10,17){\vector(1,0){4.5}}\put(14.5,17){\line(1,0){2.5}}
\put(10,24){\vector(1,0){4.5}}\put(14.5,24){\line(1,0){2.5}}
\put(13.5,13.5){\makebox(0,0){$u$}}
\put(13.5,20.5){\makebox(0,0){$u$}}
\put(8,8){\makebox(0,0){$\scr a$}}
\put(9,34.5){\makebox(0,0)[r]{$\scr m+1$}}
\put(19,8){\makebox(0,0){$\scr b$}}
\put(19,34.5){\makebox(0,0)[l]{$\scr m$}}
\put(13.5,40){\makebox(0,0)[b]{$\Phi_N^{(a,b)}(u)$}}
\end{picture}
%%%%% end subfig
}
\put(50,0){
%%%%% begin subfig 
\begin{picture}(15,35)(10,10)
\put(10,38){\vector(0,-1){4.5}}
\put(10,33.5){\vector(0,-1){7}}\put(10,26.5){\line(0,-1){16.5}}
\put(17,38){\vector(0,-1){4.5}}
\put(17,33.5){\vector(0,-1){7}}\put(17,26.5){\line(0,-1){16.5}}
\put(10,38){\vector(1,0){4.5}}\put(14.5,38){\line(1,0){2.5}}
\put(10,31){\vector(1,0){4.5}}\put(14.5,31){\line(1,0){2.5}}
\put(10,24){\vector(1,0){4.5}}\put(14.5,24){\line(1,0){2.5}}
\put(13.5,34.5){\makebox(0,0){$u$}}
\put(13.5,27.5){\makebox(0,0){$u$}}
\put(8,40){\makebox(0,0){$\scr b$}}
\put(9,13.5){\makebox(0,0)[r]{$\scr -m$}}
\put(19,40){\makebox(0,0){$\scr a$}}
\put(17.5,13.5){\makebox(0,0)[l]{$\scriptscriptstyle -m+1$}}
\put(13.5,8){\makebox(0,0)[t]{$\Phi_S^{(a,b)}(u)$}}
\end{picture}
%%%%% end subfig
}
\put(14,35){
%%%%% begin subfig 
\begin{picture}(35,15)(10,10)
\put(10,10){\vector(1,0){18.5}}
\put(28.5,10){\vector(1,0){7}}\put(35.5,10){\line(1,0){2.5}}
\put(10,17){\vector(1,0){18.5}}
\put(28.5,17){\vector(1,0){7}}\put(35.5,17){\line(1,0){2.5}}
\put(24,17){\vector(0,-1){4.5}}\put(24,12.5){\line(0,-1){2.5}}
\put(31,17){\vector(0,-1){4.5}}\put(31,12.5){\line(0,-1){2.5}}
\put(38,17){\vector(0,-1){4.5}}\put(38,12.5){\line(0,-1){2.5}}
\put(27.5,13.5){\makebox(0,0){$u$}}
\put(34.5,13.5){\makebox(0,0){$u$}}
\put(40,8){\makebox(0,0){$\scr b$}}
\put(10,7){\makebox(0,0)[l]{$\scr m$}}
\put(40,19){\makebox(0,0){$\scr a$}}
\put(10,20){\makebox(0,0)[l]{$\scr m+1$}}
\put(8,13.5){\makebox(0,0)[r]{$\Phi_W^{(b,a)}(-u)$}}
\end{picture}
%%%%% end subfig
}
\put(65,35){
%%%%% begin subfig 
\begin{picture}(35,15)(10,10)
\put(10,10){\vector(1,0){4.5}}
\put(14.5,10){\vector(1,0){7}}\put(21.5,10){\line(1,0){16.5}}
\put(10,17){\vector(1,0){4.5}}
\put(14.5,17){\vector(1,0){7}}\put(21.5,17){\line(1,0){16.5}}
\put(10,17){\vector(0,-1){4.5}}\put(10,12.5){\line(0,-1){2.5}}
\put(17,17){\vector(0,-1){4.5}}\put(17,12.5){\line(0,-1){2.5}}
\put(24,17){\vector(0,-1){4.5}}\put(24,12.5){\line(0,-1){2.5}}
\put(13.5,13.5){\makebox(0,0){$u$}}
\put(20.5,13.5){\makebox(0,0){$u$}}
\put(8,8){\makebox(0,0){$\scr a$}}
\put(38,7){\makebox(0,0)[r]{$\scr -m+1$}}
\put(8,19){\makebox(0,0){$\scr b$}}
\put(38,20){\makebox(0,0)[r]{$\scr -m$}}
\put(40,13.5){\makebox(0,0)[l]{$\Phi_E^{(b,a)}(-u)$}}
\end{picture}
%%%%% end subfig
}
\put(14,50){
%%%%% begin subfig 
\begin{picture}(35,35)(10,10)
\put(10,10){\vector(1,0){18.5}}
\put(28.5,10){\vector(1,0){7}}\put(35.5,10){\line(1,0){2.5}}
\put(24,17){\vector(1,0){4.5}}
\put(28.5,17){\vector(1,0){7}}\put(35.5,17){\line(1,0){2.5}}
\put(31,24){\vector(1,0){4.5}}\put(35.5,24){\line(1,0){2.5}}
\put(38,38){\vector(0,-1){18.5}}
\put(38,19.5){\vector(0,-1){7}}\put(38,12.5){\line(0,-1){2.5}}
\put(31,24){\vector(0,-1){4.5}}
\put(31,19.5){\vector(0,-1){7}}\put(31,12.5){\line(0,-1){2.5}}
\put(24,17){\vector(0,-1){4.5}}\put(24,12.5){\line(0,-1){2.5}}
\put(27.5,13.5){\makebox(0,0){$u$}}
\put(34.5,13.5){\makebox(0,0){$u$}}
\put(34.5,20.5){\makebox(0,0){$u$}}
\put(17,31){\makebox(0,0){$A_{NW}(u)$}}
\end{picture}
%%%%% end subfig
}
\put(14,0){
%%%%% begin subfig 
\begin{picture}(35,35)(10,10)
\put(10,38){\vector(1,0){18.5}}
\put(28.5,38){\vector(1,0){7}}\put(35.5,38){\line(1,0){2.5}}
\put(24,31){\vector(1,0){4.5}}
\put(28.5,31){\vector(1,0){7}}\put(35.5,31){\line(1,0){2.5}}
\put(31,24){\vector(1,0){4.5}}\put(35.5,24){\line(1,0){2.5}}
\put(38,38){\vector(0,-1){4.5}}
\put(38,33.5){\vector(0,-1){7}}\put(38,26.5){\line(0,-1){16.5}}
\put(31,38){\vector(0,-1){4.5}}
\put(31,33.5){\vector(0,-1){7}}\put(31,26.5){\line(0,-1){2.5}}
\put(24,38){\vector(0,-1){4.5}}\put(24,33.5){\line(0,-1){2.5}}
\put(27.5,34.5){\makebox(0,0){$u$}}
\put(34.5,34.5){\makebox(0,0){$u$}}
\put(34.5,27.5){\makebox(0,0){$u$}}
\put(17,17){\makebox(0,0){$A_{SW}(u)$}}
\end{picture}
%%%%% end subfig
}
\put(65,0){
%%%%% begin subfig 
\begin{picture}(35,35)(10,10)
\put(10,38){\vector(1,0){4.5}}
\put(14.5,38){\vector(1,0){7}}\put(21.5,38){\line(1,0){16.5}}
\put(10,31){\vector(1,0){4.5}}
\put(14.5,31){\vector(1,0){7}}\put(21.5,31){\line(1,0){2.5}}
\put(10,24){\vector(1,0){4.5}}\put(14.5,24){\line(1,0){2.5}}
\put(10,38){\vector(0,-1){4.5}}
\put(10,33.5){\vector(0,-1){7}}\put(10,26.5){\line(0,-1){16.5}}
\put(17,38){\vector(0,-1){4.5}}
\put(17,33.5){\vector(0,-1){7}}\put(17,26.5){\line(0,-1){2.5}}
\put(24,38){\vector(0,-1){4.5}}\put(24,33.5){\line(0,-1){2.5}}
\put(13.5,34.5){\makebox(0,0){$u$}}
\put(20.5,34.5){\makebox(0,0){$u$}}
\put(13.5,27.5){\makebox(0,0){$u$}}
\put(31,17){\makebox(0,0){$A_{SE}(u)$}}
\end{picture}
%%%%% end subfig
}
\put(65,50){
%%%%% begin subfig 
\begin{picture}(35,35)(10,10)
\put(10,10){\vector(1,0){4.5}}
\put(14.5,10){\vector(1,0){7}}\put(21.5,10){\line(1,0){16.5}}
\put(10,17){\vector(1,0){4.5}}
\put(14.5,17){\vector(1,0){7}}\put(21.5,17){\line(1,0){2.5}}
\put(10,24){\vector(1,0){4.5}}\put(14.5,24){\line(1,0){2.5}}
\put(10,38){\vector(0,-1){18.5}}
\put(10,19.5){\vector(0,-1){7}}\put(10,12.5){\line(0,-1){2.5}}
\put(17,24){\vector(0,-1){4.5}}
\put(17,19.5){\vector(0,-1){7}}\put(17,12.5){\line(0,-1){2.5}}
\put(24,17){\vector(0,-1){4.5}}\put(24,12.5){\line(0,-1){2.5}}
\put(13.5,13.5){\makebox(0,0){$u$}}
\put(20.5,13.5){\makebox(0,0){$u$}}
\put(13.5,20.5){\makebox(0,0){$u$}}
\put(31,31){\makebox(0,0){$A_{NE}(u)$}}
\end{picture}
%%%%% end subfig
}
\end{picture}
$$
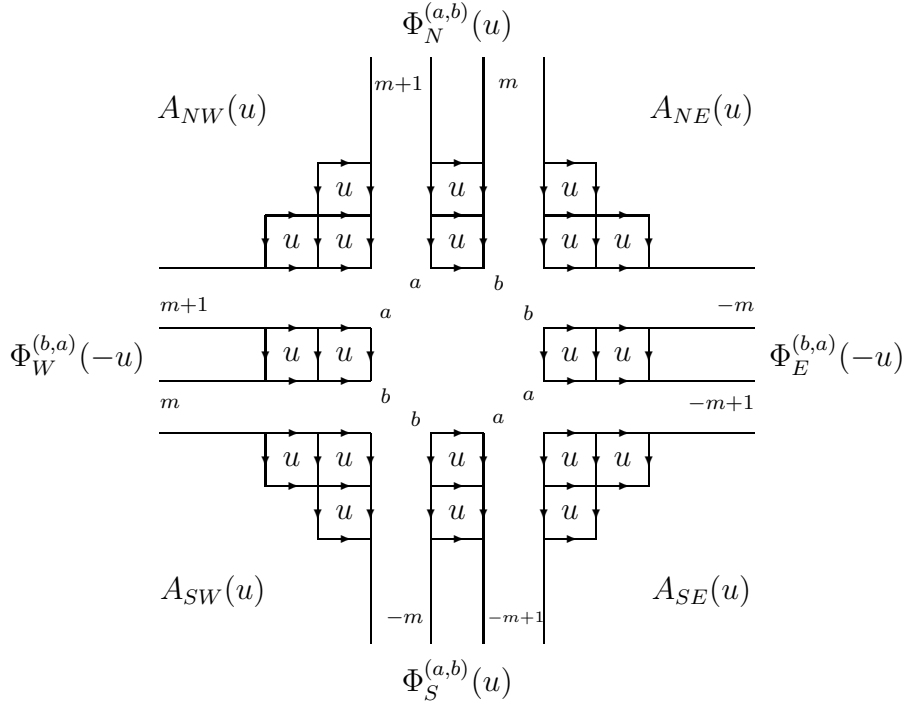
\centerline{{\bf Figure} : Corner transfer matrices and VO's of type I.}
\medskip

Formal manipulations \cite{Fo94}
show that these operators should satisfy various commutation relations. 
We have 
\be
&&\Phi_N^{(a,b)}(u)=x^{2uH^{(a)}_C}\Phi_N^{(a,b)}(0)x^{-2uH^{(b)}_C},\\
&&\Phi_W^{(a,b)}(u)=x^{2uH^{(a)}_C}\Phi_W^{(a,b)}(0)x^{-2uH^{(b)}_C}. 
\en
The VO's satisfy the commutation relations 
\bea
&& \Phi^{(a,b)}_N(u_2)  \Phi^{(b,c)}_N(u_1)
=\sum_g W\BW{a}{g}{b}{c}{u_1-u_2}
\Phi^{(a,g)}_N(u_1)\Phi^{(g,c)}_N(u_2),
\lb{VOI1}
\\
&& \Phi^{(a,b)}_W(u_2)  \Phi^{(b,c)}_W(u_1)
=\sum_g W\BW{c}{b}{g}{a}{u_1-u_2}
\Phi^{(a,g)}_W(u_1)\Phi^{(g,c)}_W(u_2),
\lb{VOI2}\\
&& \Phi^{(a,b)}_N(u_1)  \Phi^{(b,c)}_W(u_2)
=\sum_g W\BW{g}{c}{a}{b}{u_1-u_2}
\Phi^{(a,g)}_W(u_2)\Phi^{(g,c)}_N(u_1).
\lb{VOI3}
\ena
We have in addition 
\begin{eqnarray}
&&
\sum_g 
\Phi^{(a,g)}_W(u)\Phi^{(g,a)}_N(u)={\rm id},
\lb{VOI4}\\
&&
\sum_g \frac{[g]}{[a]}
\Phi^{(a,g)}_N(u) \Phi^{(g,a)}_W(k+u)={\rm id}.
\lb{VOI5}
\end{eqnarray}
The last two equations are consequences of 
the unitarity of the Boltzmann weights 
\[
\sum_{g}W\BW{a}{b}{g}{c}{u}
W\BW{a}{g}{d}{c}{-u}=\delta_{bd},
\]
and the second inversion relation
\[
\sum_g \frac{[g]}{[c]}
W\BW{d}{c}{b}{g}{-k-u}
W\BW{a}{b}{c}{g}{u}=\delta_{ad} \frac{[b]}{[d]}.
\]

In Section 3 we shall present a bosonic realization of the relations 
\eqref{VOI1}--\eqref{VOI5}.

%%%%%%%%%%%%%%%%%%%%%%%%%%%%
%  2.3                     %
%%%%%%%%%%%%%%%%%%%%%%%%%%%%
\subsection{Local height probabilities}\lb{subsec:2.3}

Consider successive $n$ sites $i_1,\cdots,i_n$ on a row of the lattice, 
numbered from right to left.  
We regard $i_1$ as the reference site used to label the ground states. 
Fixing the boundary heights to the ground state $m$,  
we let $P_{a_n,\cdots,a_1}(m)$ denote the joint 
probability that the height variable $l_{i_j}$ takes 
the value $a_j$, $j=1,\cdots,n$.  
These are the $n$ point local height probabilities (LHP's). 
In terms of the CTM and VO's, they can be expressed as 
\be
&&P_{a_n,\cdots,a_1}(m)=\frac{1}{Z_m}[a_1]\\
&&~~\times
{\rm tr}_{\H_{m,a_1}}\Bigl(x^{2kH^{(a_1)}_C}
\Phi^{(a_1,a_{2})}_W(0)\cdots \Phi^{(a_{n-1},a_n)}_W(0)
\Phi^{(a_n,a_{n-1})}_N(0)\cdots \Phi^{(a_2,a_1)}_N(0)
\Bigr), \\
&&Z_m=\sum_{a=1}^{k+1}[a]\,
{\rm tr}_{\H_{m,a}}\bigl(x^{2kH^{(a)}_C}\bigr).
\en
It turns out that (with an appropriate normalization of $H^{(a)}_C$)
$Z=Z_m$ is independent of $m$. 
Obviously $P_{a_n,\cdots,a_1}(m)=0$ for $a_1\equiv m \bmod 2$.

The following result for $n=1$ is 
due to Andrews, Baxter and Forrester \cite{ABF}.
\bea
P_a(m)=\frac{[a]\,{\rm tr}_{\H_{m,a}}\bigl(x^{2kH^{(a)}_C}\bigr)}{Z}
=
x^{\frac{k+2}{4}}[a]\,
c^{\Lambda_k(a-1)}_{\Lambda_k(m)}(\tau),  
\lb{LHP1}
\ena
where $c^{\Lambda_k(l)}_{\Lambda_k(m)}(\tau)$  
stands for the string function \cite{KaPe} for integrable $\slth$-modules.
Explicitly it is  
\be 
&&c^{\Lambda_k(l)}_{\Lambda_k(m)}(\tau) 
=\eta(\tau)^{-3}
\\
&&\qquad \times
\Bigl(\sum_{n_1\geq |n_2|}-\sum_{-n_1> |n_2|}\Bigr)
(-1)^{2n_1}
\bigl(x^{2k}\bigr)^{\frac{(l+1+2n_1(k+2))^2}{4(k+2)}-\frac{(m+2n_2k)^2}{4k}}
\en
for $l\equiv m\bmod 2$ and 
$c^{\Lambda_k(l)}_{\Lambda_k(m)}(\tau)=0$ for $l\not\equiv m\bmod 2$.
In the above, we set $x^{2k}=e^{2\pi i \tau}$, 
$\eta(\tau)=(x^{2k})^{1/24}(x^{2k};x^{2k})_\infty$,    
$\Lambda_k(j)=(k-j)\Lambda_0+ j \Lambda_1$, 
and the sum is taken over $n_1,n_2\in\frac{1}{2}\Z$ with
 $n_1-n_2\in \Z$. 
We have
\be
&&P_a(m)=0\quad  (a\equiv m\bmod 2),
\\
&&\sum_{a=1}^{k+1}P_a(m)=1,  
\\
&&P_a(-m)=P_a(m).  
\en

As observed in \cite{LP96} for regime III, 
the nearest neighbor LHP's ($n=2$) 
can be written in terms of the one point LHP. 
To see this it suffices to note  
the following relations which are obvious consequences of the definition.  
\be
&&P_{1,2}(m)=P_1(m+1),\qquad
P_{2,1}(m)=P_1(m),\\
&&P_{a-1,a}(m)+P_{a+1,a}(m)=P_a(m), \\ 
&&P_{a,a-1}(m)+P_{a,a+1}(m)=P_a(m+1). 
\en
These properties fix $P_{a\pm 1,a}(m)$ uniquely in terms of the 
one point LHP \eqref{LHP1} as follows. 
\bea
&&P_{a+1,a}(m)=\sum_{1\le s\le a \atop s\equiv a\bmod 2}P_s(m)
-\sum_{1\le s\le a \atop s\not\equiv a\bmod 2}P_s(m+1), 
\lb{LHP2}\\
&&P_{a-1,a}(m)=-\sum_{1\le s< a \atop s\equiv a\bmod 2}P_s(m)
+\sum_{1\le s < a \atop s\not\equiv a\bmod 2}P_s(m+1).
\lb{LHP3}
\ena
The right hand sides of \eqref{LHP2}--\eqref{LHP3} 
automatically satisfy 
 $P_{k+2,k+1}(m)=P_{k+1,k+2}(m)=0$, 
as it should be. 

\setcounter{section}{2}
\setcounter{equation}{0}
%%%%%%%%%%%%%%%%%%%%%%%%%%%%%%%%%%%%%%%%%%%%%%%%%%%%%%%%%%%%%%%
%                                                             %
%  3. Bosonization                                            %
%                                                             %
%%%%%%%%%%%%%%%%%%%%%%%%%%%%%%%%%%%%%%%%%%%%%%%%%%%%%%%%%%%%%%%
\section{Free Field Realization} \lb{sec:3}

The $q$-deformation of the Wakimoto module over the affine Lie algebra
$\slth$ was found by several authors \cite{Ma94}
(See refs. in \cite{Ma94} as to the other variants of 
free field realizations.), 
using three kinds of bosonic fields.
The deformed parafermion theory is obtained by
dropping one of these fields.
In this section we give some details of the latter theory,
and interpret it as a free field realization of the ABF model in regime II.
A part of the results are given in \cite{Konno}. 
The various relations 
given in Propositions \ref{prop:3.1}, 
\ref{prop:3.3}--\ref{prop:3.4} are new. 

%%%%%%%%%%%%%%%%%%%%%%%%%%%%
%  3.1                     %
%%%%%%%%%%%%%%%%%%%%%%%%%%%%
\subsection{Basic operators}\lb{subsec:3.1}

The main objects in the deformed parafermion theory are
the vertex operators (VO's) of type I
\bea
\Phi^*_\pm(u), \quad \Phi_\pm(u),
\lb{I}
\ena
and of type II
\bea
\Psi_+(u)=\Psi^\dagger(u),\quad \Psi_-(u)=\Psi(u).
\lb{II}
\ena
The VO's of type I \eqref{I} are a lattice counterpart of the
simplest primary fields in conformal theory,
while those of type II \eqref{II} corresponds to the parafermions.
In addition, there are also some auxiliary operators: the so-called
screening current and the $\xi$-$\eta$ system which we denote by
\be
S(u),\quad \xi(u),\quad \eta(u),
\en
respectively.
All these operators act on the direct sum of Fock modules
$\F=\oplus_{m\equiv l\bmod 2 }\F_{m,l}$ labeled by $m,l\in \Z$.
Their explicit formulas are given in Appendix \ref{sec:app1}.
Here we only mention their basic features.

The operators
\be
\Psi_{\pm}(u),~\Phi^*_-(u),~\Phi_-(u),~S(u),~\xi(u),~\eta(u)
\en
are either an exponential of a bosonic field or a sum of two such terms.
In contrast, the formulas for $\Phi^*_+(u),~\Phi_+(u)$ comprise integrals.
To elaborate on the last point, let us introduce the screening charge
\be
X(u)=\oint_{C_X(z)}\underline{ dz}'
S(u')\frac{[u-u'-\frac{1}{2}-P_1]}{[u-u'-\frac{1}{2}]}\,.
%\lb{X(u)}
\en
Here we set $z=x^{2u}$, $z'=x^{2u'}$, $\underline{dz}'=dz'/(2 \pi i z')$,
and $P_1$ denotes a `zero-mode' operator.
For more details we refer to Appendix \ref{sec:app1}. 
The contour $C_X(z)$ is a simple closed curve encircling
the origin counterclockwise in the $z'$-plane, such that
$z'=zx^{-1+ 2(k+2) n}$ ($n\geq 0$) are inside $C_X(z)$ and
$z'=zx^{-1+ 2(k+2) n}$ ($n< 0$) are outside $C_X(z)$.
Using the screening charge, the `plus'
components of the type I VO's are given by
\be
&&\Phi_{+}(u)=\Phi_{-}(u)X(u+k+1),
\\
&& \Phi^*_{+}(u)=-\Phi^*_{-}(u) X(u),
\en
where the product $A(u)B(v)$ is defined to be the analytic continuation
from the domain $|x^{2u}|\gg|x^{2v}|$.

Each of the operators
$\Phi_{\pm}(u)$, $\Phi^*_{\pm}(u)$, $\Psi_{\pm}(u)$,
$S(u)$, $\xi(u)$, $\eta(u)$
sends one Fock space $\F_{m,l}$ to
another $\F_{m',l'}$ with a different label $(m',l')$.
The change of label is listed in the following table.
\medskip

$$
\begin{tabular}{|c|c|c|c|c|c|c|}
\hline
 &\mathstrut $\Phi_{\pm}(u)$ &$\Phi^*_{\pm}(u)$ &$\Psi_{\pm}(u)$ &$S(u)$
& $\xi(u)$ & $\eta(u)$ \\
\hline
$\F_{m',l'}$&$\F_{m+1,l\mp1}$&$\F_{m-1,l\mp1}$&$\F_{m\mp2,l}$&
$\F_{m,l-2}$
&$\F_{m-k,l-(k+2)}$&$\F_{m+k,l+k+2}$\\
\hline
\end{tabular}
$$
\medskip

%%%%%%%%%%%%%%%%%%%%%%%%%%%%
%  3.2                     %
%%%%%%%%%%%%%%%%%%%%%%%%%%%%
\subsection{Commutation relations}\lb{subsec:3.2}

In the next three propositions
we state the commutation relations among the operators
 $\Phi_{\pm}(u)$, $\Phi^*_{\pm}(u)$ and $\Psi_\pm(u)$.
Set
\be
P=P_1+1.
\en

\begin{prop}\lb{prop:3.1}
The operators $\Phi_{\pm}(u)$, $\Phi^*_{\pm}(u)$
satisfy the commutation relations
\be
&&
\Phi_{\ve_2}(u_2)\Phi_{\ve_1}(u_1)
=\sum
W\BW{P}{P+\ve_1'}{P+\ve_2}{P+\ve_1+\ve_2}{u_1-u_2}
\Phi_{\ve_1'}(u_1)\Phi_{\ve_2'}(u_2),
\\
&&
\Phi^{*}_{\ve_2}(u_2)  \Phi^{*}_{\ve_1}(u_1)
=\sum
W\BW{P+\ve_1+\ve_2}{P+\ve_2}{P+\ve_1'}{P}{u_1-u_2}
\Phi^{*}_{\ve_1'}(u_1)\Phi^{*}_{\ve_2'}(u_2),
\\
&&
\Phi_{\ve_1}(u_1)\Phi^{*}_{\ve_2}(u_2)
=\sum
W\BW{P+\ve_2'}{P+\ve_1+\ve_2}{P}{P+\ve_1}{u_1-u_2}
\Phi^{*}_{\ve_2'}(u_2)\Phi_{\ve_1'}(u_1),
\en
where the Boltzmann weights are given in \eqref{BW1}--\eqref{BW3}.
The sums are taken over $\ve_1',\ve_2'=\pm 1$ such that
$\ve_1'+\ve_2'=\ve_1+\ve_2$.
Moreover we have the inversion relations
\bea
&&
\sum_{\ve}
\Phi^{*}_\ve(u)\Phi_{-\ve}(u)={\rm id},
\lb{finv}\\
&&
\sum_{\ve}
\Phi_\ve(u)[P] \Phi^{*}_{-\ve}(k+u)=[P]\,{\rm id}.
\lb{sinv}
\ena
\end{prop}

In the next propositions we set
\be
[u]^*=x^{\frac{u^2}{k}-u}\Theta_{x^{2k}}(x^{2u}).
\en
\begin{prop}\lb{prop:3.2}
The following commutation relations hold.
\begin{eqnarray*}
&& \Psi_\pm (u_2) \Psi_\pm (u_1) =
\frac{[u_1-u_2+1]^*}{[u_1-u_2-1]^* }
 \Psi_\pm (u_1) \Psi_\pm (u_2),\\
&& \Psi_\pm (u_2) \Psi_\mp (u_1) =
\frac{[u_1-u_2-1+\frac{k}{2}]^*}{[u_1-u_2+1+\frac{k}{2}]^* }
\Psi_\mp (u_1) \Psi_\pm (u_2).
\end{eqnarray*}
Moreover, as 
$u_2\rightarrow u_1-\frac{k}{2}$ we have
\be
\Psi_{\pm}(u_2)\Psi_{\mp}(u_1)
=
\frac{k}{2\pi}\frac{g^*}{u_2-u_1+\frac{k}{2}}\times\id +O(1),
\en
where 
$\displaystyle g^*=\frac{\pi}{k\log x}\frac{(x^{2k-2};x^{2k})_\infty}
{(x^{2};x^{2k})_\infty}$.
\end{prop}

\begin{prop}\lb{prop:3.3}
\begin{eqnarray*}
&& \Psi_+ (u_2) \Phi_{\pm} (u_1) =
\frac{[u_1-u_2-\frac{1}{2}+\frac{k}{2}]^*}
{[u_1-u_2+\frac{1}{2}+\frac{k}{2}]^* }
 \Phi_{\pm} (u_1) \Psi_+ (u_2),
\\
&& \Psi_+ (u_2) \Phi^*_{\pm} (u_1) =
\frac{[u_1-u_2+\frac{1}{2}+\frac{k}{2}]^*}
{[u_1-u_2-\frac{1}{2}+\frac{k}{2}]^*}
 \Phi^*_{\pm} (u_1) \Psi_+ (u_2),
\\
&& \Psi_- (u_2) \Phi_{\pm} (u_1) =
\frac{[u_1-u_2+\frac{1}{2}]^*}{[u_1-u_2-\frac{1}{2}]^* }
 \Phi_{\pm} (u_1) \Psi_- (u_2),
\\
&&
\Psi_- (u_2) \Phi^*_{\pm} (u_1) =
\frac{[u_1-u_2-\frac{1}{2}]^*}{[u_1-u_2+\frac{1}{2}]^* }
 \Phi^*_{\pm} (u_1) \Psi_- (u_2).
\end{eqnarray*}
\end{prop}

%%%%%%%%%%%%%%%%%%%%%%%%%%%%
%  3.3                     %
%%%%%%%%%%%%%%%%%%%%%%%%%%%%
\subsection{Free fields resolution}\lb{subsec:3.3}

The space of states for the deformed parafermion theory 
is constructed following a procedure 
well known in conformal field theory \cite{BeFe90}.
This is done in two steps.
The first step is to introduce a certain subspace $\Ft_{m,l}$ of the Fock
space $\F_{m,l}$ using $\eta(u)$.
The second step is to consider a complex consisting of
these $\Ft_{m,l}$, in which the coboundary maps are
given by powers of the screening charge $X(u)$.
The space of states is then defined as the $0$-th cohomology 
of this complex.

Let
\be
\eta_0:\F_{m,l}\longrightarrow \F_{m+k,l+k+2}
\en
denote the zeroth Fourier coefficient of $\eta(u)$.
It is well defined provided $m\equiv l\bmod 2$.
We set
\be
\Ft_{m,l}={\rm Ker }~\eta_0\bigl|_{\F_{m,l}}.
\en
Then we have the following resolution of $\Ft_{m,l}$
by Fock spaces (see Appendix \ref{sec:app2})
\bea
&&\lb{eta0}\\
&&0\longrightarrow \Ft_{m,l}\longrightarrow
\F_{m,l}\overset{\eta_0}{\longrightarrow }
\F_{m+k,l+k+2} \overset{\eta_0}{\longrightarrow }
\F_{m+2k,l+2(k+2)}\overset{\eta_0}{\longrightarrow }
\cdots,
\nn\\
&&\lb{eta00}\\
&&
\cdots \overset{\eta_0}{\longrightarrow }\F_{m-2k,l-2(k+2)}
 \overset{\eta_0}{\longrightarrow }\F_{m-k,l-(k+2)}
\overset{\eta_0} \longrightarrow \Ft_{m,l}\longrightarrow 0.
\nn
\ena

The following proposition shows that the operators
$S(u)$, $\Psi_\pm(u)$, $\Phi_{\pm}(u)$ and
$\Phi^*_{\pm}(u)$ have a well-defined action on the subspace $\Ft_{m,l}$.
\begin{prop} We have
\be
&&\eta_0 S(u)=-S(u) \eta_0,\\
&&\eta_0 \Psi_{\pm}(u)=- \Psi_{\pm}(u)\eta_0, \\
&&\eta_0 \Phi_{\pm}(u)=  \Phi_{\pm}(u)\eta_0, \\
&&\eta_0 \Phi^*_{\pm}(u)= \Phi^*_{\pm}(u)\eta_0.
\en
\end{prop}

Next we fix $l,m\in\Z$ with $0\le l\le k$, $l\equiv m\bmod 2$.
Consider a sequence $\Cc_{m,l}$
\bea
&&\lb{resol}\\
&&
\cdots ~{\buildrel X_{-2} \over \longrightarrow }~
\Ft_{m,-l-2+2(k+2)} ~{\buildrel X_{-1} \over \longrightarrow }~
\Ft_{m,l} ~{\buildrel X_0 \over \longrightarrow }~
\Ft_{m,-l-2} ~{\buildrel X_1 \over \longrightarrow }~
\Ft_{m,l-2(k+2)} ~{\buildrel X_2 \over \longrightarrow }~
\cdots,
\nn
\ena
defined by appropriate powers of the screening charge $X(u)$, i.e.,
\be
&&X_{2j}=X(u)^{l+1}~:~\Ft_{m,l-2j(k+2)}\longrightarrow
\Ft_{m,-l-2-2j(k+2)}, \\
&&X_{2j+1}=X(u)^{k-l+1}~:~\Ft_{m,-l-2-2j(k+2)}\longrightarrow
\Ft_{m,l-2(j+1)(k+2)}.
\en
The following can be shown in exactly the same way as in \cite{JLMP}.
\begin{prop} The maps $X_j$ are independent of $u$, and
 $\Cc_{m,l}$ is a cochain complex:
\be
&&X_jX_{j-1}=0 \qquad (j\in\Z).
\en
\end{prop}
The following statement concerning the cohomology
group of this complex seems quite plausible (cf. \cite{BeFe90}).
\bea
H^j(\Cc_{m,l})=0 \qquad (j\neq 0).
\lb{coho}
\ena
As we do not have a rigorous mathematical proof,
we assume henceforth the validity of \eqref{coho}.
By the Euler-Poincar\'e principle,
the character of the remaining cohomology $H^0(\Cc_{m,l})$ then becomes
\be
&&{\rm tr}_{H^0(\Cc_{m,l})}\bigl(x^{2kD}\bigr)=
x^{\frac{k^2}{4(k+2)}}(x^{2k};x^{2k})_\infty\,
c^{\Lambda_k(l)}_{\Lambda_k(m)}(\tau),
\\
&&Z=\sum_{0\le l\le k\atop l \equiv m \bmod 2}
[l+1]\,{\rm tr}_{H^0(\Cc_{m,l})}\bigl(x^{2kD}\bigr)=
x^{-\frac{k+1}{k+2}}(x^{2k};x^{2k})_\infty,
\en
giving the same formula as that of the one point LHP \eqref{LHP1}.

\begin{prop}\lb{prop:3.4}
\be
&&X_j \Psi_\pm(u)=\Psi_\pm(u)X_j,
\\
&&X_j \Phi_{\pm}(u)=\Phi_{\mp}(u)X_j,
\\
&&X_j \Phi^*_{\pm}(u)=\Phi^*_{\mp}(u)X_j.
\en
\end{prop}
This proposition ensures that
 $\Phi_{\pm}(u)$, $\Phi^*_{\pm}(u)$ and $\Psi_\pm(u)$ give rise to
well-defined operators on the cohomology $H^0(\Cc_{m,l})$.
We abuse the notation and denote them by the same letters.

%%%%%%%%%%%%%%%%%%%%%%%%%%%%
%  3.4                     %
%%%%%%%%%%%%%%%%%%%%%%%%%%%%
\subsection{Identification with lattice theory}

The construction of this Section is related to the
lattice theory in Section 2 as follows.
We make the following identification:
\begin{enumerate}
\item The space of states of the lattice model
(with central height $a$ and boundary condition $m$)
with the $0$-th cohomology of $\Cc_{m,a-1}$,
\be
\H_{m,a}=H^0(\Cc_{m,a-1}).
\en
\item The corner Hamiltonian $H^{(a)}_C$
with the grading operator $D$.
\item The half-infinite transfer matrices with the
type I VO's
\be
&&\Phi^{(a-\varepsilon,a)}_N(u)=\Phi_\varepsilon(u)\bigl|_{H^0(\Cc_{m,a-1})},\\
&&\Phi^{(a-\varepsilon,a)}_W(u)=\Phi^*_\varepsilon(u)\bigl|_{H^0(\Cc_{m,a-1})}.
\en
\item The creation/annihilation operators of particles
and anti-particles with the type II VO's $\Psi_\pm(u)$.
\end{enumerate}
As we already mentioned,
with this identification the characters of the two spaces match.
The commutation relations for VO's of type I expected from the
lattice theory \eqref{VOI1}--\eqref{VOI5}
are recovered in Proposition \ref{prop:3.1}.
In Section \ref{sec:5}, we will comment about the agreement 
with the known results on the excitation spectrum  
and the $S$-matrix \cite{BaRe90} in the scaling limit.

Correlation functions and form factors of local operators are given
as traces of operators acting on the physical space $H^0(\Cc_{m,l})$.
Let $\mathcal{O}$ stand for $x^{2kD}$ times a product of operators of the form
$\Psi_{\varepsilon}(u)$,
$\Phi_\varepsilon(u)$ and $\Phi^*_\varepsilon(u)$.
For simplicity we consider the case where there are an equal number of
$\Psi_{+}(u)$'s and  $\Psi_{-}(u)$'s.
Proposition \ref{prop:3.4} shows that
\be
&&\mathcal{O}_i\eta_0=\eta_0\mathcal{O}_i,
\\
&&X_j \mathcal{O}_i=\mathcal{O}_{1-i} X_j
\qquad (j\equiv i \bmod 2),
\en
where
$\mathcal{O}_0=\mathcal{O}$, and
$\mathcal{O}_1$ signifies the operator
obtained from $\mathcal{O}$ by negating
the indices $\varepsilon$
of $\Phi_\varepsilon(u)$ and $\Phi^*_\varepsilon(u)$.
The resolutions \eqref{eta0},\eqref{eta00},\eqref{resol} 
afford a procedure for computing the trace as follows.
\bea
&&\tr_{H^0(\Cc_{m,l})}(\mathcal{O})
=\sum_{s\in\Z}
\tr_{\Ft_{m,l-2s(k+2)}}(\mathcal{O}_0)
-
\sum_{s\in\Z}
\tr_{\Ft_{m,-l-2-2s(k+2)}}(\mathcal{O}_1),
\lb{resol1}
\ena
\bea
\tr_{\Ft_{m,l}}(\mathcal{O}_i)
&\!\!=\!\!&\sum_{n\ge 0}(-1)^n\, \tr_{\F_{m+kn,l+(k+2)n}}(\mathcal{O}_i)
\lb{resol2}\\
&\!\!=\!\!&
-\sum_{n< 0}(-1)^n \,\tr_{\F_{m+kn,l+(k+2)n}}(\mathcal{O}_i).
\nn
\ena
Taking 
$\mathcal{O}=x^{2kD}\Phi^*_{\varepsilon}(u)\Phi_{-\varepsilon}(u)$,
we have verified directly that the formula thus obtained reproduces 
the LHP \eqref{LHP2}--\eqref{LHP3} for neighboring lattice sites. 

%%%%%%%%%%%%%%%%%%%%%%%%%%%%
%  3.5                     %
%%%%%%%%%%%%%%%%%%%%%%%%%%%%
\subsection{Fusion of parafermions and $W$ currents} 

The type II VO's $\Psi_-(u)$, $\Psi_+(u)$ 
play the role of creation operators of excitations over the ground state. 
Besides them, there are altogether $k-1$ kinds of `elementary' particles 
in regime II, as expected from the level-rank duality \cite{JMO1}.  
The standard fusion procedure provides us with a free field realization 
$\Psi_a(u)$ defined recursively as 
\be
  &&\Psi_1(u'-\frac{a}{2})\Psi_a(u+\frac{1}{2})\\
  &=\!\!&
  \frac{1}{1-\frac{z}{z'}}x^{\frac{a(a+1)}{k}}
  \frac{(x^{2a},x^{2a+2k+2};x^{2k})_{\infty}}
       {(x^{2k},x^{2k+2};x^{2k})_{\infty}}\,
  \Psi_{a+1}(u)+O(1)\quad (u'\rightarrow u), 
\en
where $\Psi_1(u)=-(x-x^{-1})\Psi_-(u)$. 
A more explicit expression is given in \eqref{fusP}. 
We have the following commutation relations
\bea
  &&\Psi_a(u_2)\Phi_\pm(u_1)=
  \frac{[u_1-u_2+{a\over 2}]^*}{[u_1-u_2-{a\over 2}]^*}
  \Phi_\pm(u_1)\Psi_a(u_2),\nn\\
  &&\Psi_a(u_2)\Psi_b(u_1)
  \lb{PsiaPsib}\\
  &=\!\!&
  \prod_{s=1}^a\prod_{s'=1}^b
  \frac{[u_1-u_2+1+\frac{a-b}{2}-(s-s')]^*}{[u_1-u_2-1+\frac{a-b}{2}-(s-s')]^*}
  \cdot\Psi_b(u_1)\Psi_a(u_2)\nn\\
  &=\!\!&
  \frac{[u+\frac{a+b}{2}]^*}{[u-\frac{a+b}{2}]^*}
  \frac{[u+\frac{|a-b|}{2}]^*}{[u-\frac{|a-b|}{2}]^*}
  \prod_{s=1}^{\min(a,b)-1}
  \frac{{[u+\frac{|a-b|}{2}+s]^*}^2}{{[u-\frac{|a-b|}{2}-s]^*}^2}
  \cdot\Psi_b(u_1)\Psi_a(u_2),\nn
\ena
where $u=u_1-u_2$.
These relations agree with the known results about the 
excitation spectra \eqref{DiRe} 
and the scattering matrices \eqref{Sab} in the scaling field theory 
discussed in Section \ref{sec:5}. 

The operator $\Psi_k(u)$ commutes with $\Psi_{\pm}(v)$ and $X_j$, 
commutes or anticommutes with $\eta_0$, 
and anticommutes with $\Phi_{\pm}(v)$. 
It can be shown that it is also invertible. 
We expect that  on the cohomology it is independent of $u$ and 
defines an isomorphism $\iota:H^0(\Cc_{m,l}) \simeq H^0(\Cc_{m+2k,l})$.  
We also expect that $\Psi_{k-1}(u)$ 
defines the same operator as $\iota\circ \Psi_+(u)$.   

The level-rank duality also suggests the existence of 
the deformed $W_{k}$ currents \cite{AKOS96,FeFr95} 
in the parafermionic description of the ABF models. 
Indeed, the first deformed $W$ current $W^1(u)$ 
can be obtained either as 
fusion of type I VO's ($z_2\rightarrow x^{2(k+2)}z_1$) 
\be
 &&\Phi_-(u_2)\Phi^*_+(u_1)\\
  &=\!\!&\Bigl(1-x^{2(k+2)}\frac{z_1}{z_2}\Bigr)
  \frac{x^{\frac{1}{k}}}{[k+1]_x}
  \frac{(x^{2k},x^{-4};x^{2k})_{\infty}}{(x^{-2},x^{-2};x^{2k})_{\infty}}\,
  W^1(u_1)+\cdots,
\en
or as that of type II VO's ($z_2\rightarrow x^{k+2}z_1$)
\be
  &&\Psi_-(u_2)\Psi_+(u_1)\\
  &=\!\!&\frac{1}{1-x^{k+2}\frac{z_1}{z_2}}
  \frac{x^{\frac{2}{k}+1}}{(x-x^{-1})^2[k+1]_x}
  \frac{(x^{-4};x^{2k})_{\infty}}{(x^{2k};x^{2k})_{\infty}}\,
  W^1(u_1-\sfrac{k+1}{2})+\cdots.
\en
Explicit formulas of 
$W^1(u)$, 
as well as those of the higher currents $W^j(u)$ ($j=2,3,\cdots$),  
are given in Appendix \ref{sec:app1}. 
We expect that these $W^j(u)$ ($j\geq 1$) generate the same 
deformed $W$ algebra for $\slN$ in \cite{AKOS96}, 
under the following identification
\be
  N=k,\quad q=x^{2(k+1)},\quad t=x^{2(k+2)}\quad (p=qt^{-1}=x^{-2}).
\en
We have checked 
a part of the relations for $W^j(u)$ (eq.(8) in \cite{AKOS96}), 
but have not verified such relations 
as $W^k(u)=1$ which are expected to hold only at the level
of the cohomology.

\setcounter{section}{3}
\setcounter{equation}{0}
%%%%%%%%%%%%%%%%%%%%%%%%%%%%%%%%%%%%%%%%%%%%%%%%%%%%%%%%%%%%%%%
%                                                             %
%  4. Form Factors                                            %
%                                                             %
%%%%%%%%%%%%%%%%%%%%%%%%%%%%%%%%%%%%%%%%%%%%%%%%%%%%%%%%%%%%%%%
\section{Form Factors} \lb{sec:4}

%%%%%%%%%%%%%%%%%%%%%%%%%%%%
%  4.1                     %
%%%%%%%%%%%%%%%%%%%%%%%%%%%%
\subsection{Traces of type II operators}\lb{subsec:4.1}

As a simplest example of form factors of local operators,
let us consider the quantity
\bea
&&
Q_{a}^{(n,n)}(m)=\bigl(Zg^{*n}\bigr)^{-1}[a]
\lb{Qbar}\\
&&\quad \times
\tr_{\H_{m,a}}\Bigl(x^{2kD}\Psi_{+}(v_1)\cdots \Psi_{+}(v_n)
\Psi_{-}(v'_{1})\cdots\Psi_{-}(v'_{n})\Bigr).
\nn
\ena
To simplify the notation, we have not exhibited explicitly 
the dependence of $Q_{a}^{(n,n)}(m)$ on the parameters 
$v_j,v_j'$. 
For the same reason we will often suppress the superscript $(n,n)$.
Note that $Q_a(m)=0$ for $a\equiv m\bmod 2$.

Hereafter we set 
\be
\tau=-\frac{ik}{\pi}\log x\,.
\en
It is a standard task to compute traces of bosonic operators over 
the Fock space. 
After the working outlined in Appendix \ref{sec:app3}, 
we find the following expression.
\bea
  &&Q_{a}(m)=
  \sum_{\mu_1,\cdots,\mu_n=\pm 1 \atop \nu_1,\cdots,\nu_n=\pm 1}
  \kappa(v,\mu,\nu)
  \prod_{i=1}^{n} \mu_i\nu_i
  \lb{Qlat}\\
  &&\quad\times\prod_{1 \leq i<j\leq n}
  \frac{[v_j-v_i+\frac{\mu_i-\mu_j}{2}]^*}{[v_j-v_i-1]^*}
  \frac{F(v_{j}-v_i)}{F(0)}
  \nn\\
  &&\quad\times\prod_{1 \leq i<j\leq n}
  \frac{[v'_{j}-v'_{i}-\frac{\nu_i-\nu_j}{2}]^*}{[v'_{j}-v'_{i}-1]^*}
  \frac{F(v'_{j}-v'_{i})}{F(0)}
  \nn\\
  &&\quad\times\prod_{i,j=1}^{n}
  \frac{[v'_{j}-v_{i}-\frac{k}{2}-\frac{\mu_i+\nu_j}{2}]^*}
  {[v'_{j}-v_{i}-\frac{k}{2}]^*}
  \frac{F(0)}{F(v'_{j}-v_i-\frac{k}{2})}
  \nn\\
  &&\quad\times\,\Gamma_{m,a-1}\Bigl(\frac{\tau}{2k}(\mu-\nu),
  \frac{\tau}{k}\bigl(\frac{2v}{k}-\frac{\mu+\nu}{2}\bigr)\Bigl|\tau\Bigr).
  \nn
\ena
In the above, we have set
$v=\sum_{i=1}^n(v'_i-v_i)$,
$\mu=\sum_{i=1}^n\mu_i$, $\nu=\sum_{i=1}^n\nu_i$.
The function 
$\Gamma_{m,l}(y_1,y_2|\tau)=\Gamma^{(0,0)}_{m,l}(y_1,y_2|\tau)$ 
is defined in \eqref{Gam}. 
The functions $F(v)$ and 
$\kappa(v,\mu,\nu)$ are defined as follows.
\bea
&&
F(v)=
\frac{(x^{2(k+1+v)},x^{2(k+1-v)};x^{2k},x^{2k})_\infty}
{(x^{2(k-1+v)},x^{2(k-1-v)};x^{2k},x^{2k})_\infty},
\lb{F(w)}
\\
&&
\kappa(v,\mu,\nu)=
\frac{[a]}{Z}
\Bigl(-i\tau \frac{\eta(\tau)^3}{[1]^{*}}\Bigr)^n
x^{(\frac{\mu+\nu}{k}-\frac{2}{k}n)v+
\frac{n^2}{2}+\frac{\mu^2+\nu^2}{4k}-
\frac{kn}{4}-\frac{k+1}{2k}\mu\nu+\frac{kc}{12}}
\nn
\ena
with $c=2(k-1)/(k+2)$. 

Using \eqref{Gamsym1},\eqref{Gamsym2} 
we see that $Q_a(m+2k)=Q_{k+2-a}(m+k)=Q_a(m)$. 
\medskip

\noindent{\it Remark.}\quad
We see from \eqref{contr1}, \eqref{contr2} that 
the contribution of the oscillator part to the trace \eqref{Qbar} 
is convergent if 
$x^{\varepsilon_2-\varepsilon_1}<|z_1/z_2|
<x^{\varepsilon_2-\varepsilon_1-2k}$. 
For $k>2$ there is a non-empty domain of convergence
common to all $\varepsilon_1, \varepsilon_2$.   
The case $k=2$ of the Ising model is exceptional and needs
a separate treatment by analytic continuation. 
{\it For the rest of the paper we assume $k\ge 3$.} 

%%%%%%%%%%%%%%%%%%%%%%%%%%%%
%  4.2                     %
%%%%%%%%%%%%%%%%%%%%%%%%%%%%
\subsection{Neighboring heights}\lb{subsec:4.2}

As the next example, let us consider the trace involving two type I VO's,
\be
Q_{b,a}(m)&\!\!=\!\!&\bigl(Zg^{*n}\bigr)^{-1}[a]\,
\tr_{\H_{m,a}}\Bigl(x^{2kD}\Phi^{*(a,b)}(u)\Phi^{(b,a)}(u)
\\
&&\qquad\times \Psi_{+}(v_1)\cdots \Psi_{+}(v_n)
\Psi_{-}(v'_{n})\cdots \Psi_{-}(v'_{n})\Bigr).
\en
Here
$\Phi^{(a-\varepsilon,a)}(u)$ means $\Phi_{\varepsilon}(u)|_{\H_{m,a}}$,
and similarly for $\Phi^{*(a,b)}(u)$.

We have already mentioned that the two point LHP for the neighboring height
variables \eqref{LHP2},\eqref{LHP3}
can be expressed simply in terms of one point LHP's.
Let us apply the same argument to compute $Q_{b,a}(m)$.
The first inversion identity \eqref{finv} entails that
\bea
Q_{a+1,a}(m)+Q_{a-1,a}(m)=Q_{a}(m),
\lb{Q1}
\ena
with $Q_{a}(m)$ given by \eqref{Qbar}.
On the other hand, the second inversion relation \eqref{sinv}
together with the cyclicity of the trace implies
\bea
Q_{b,b-1}(m)+Q_{b,b+1}(m)=Q_b(m+1)\,G.
\lb{Q2}
\ena
The factor
\be
&&
G=\prod_{j=1}^n \frac{[u-v_j+\frac12+\frac{k}{2}]^*}
  {[u-v_j-\frac12+\frac{k}{2}]^*}
\frac{[u-v'_{j}-\frac{1}{2}]^*}{[u-v'_{j}+\frac{1}{2}]^*}
\en
arises from the commutation relations between
type I and type II VO's.
Solving the relations \eqref{Q1},\eqref{Q2} under the
condition $Q_{2,1}(m)=Q_1(m)$, $Q_{1,2}(m)=Q_1(m+1)$
(which can be verified using the integral representations),
we obtain
\be
Q_{a+1,a}(m)=
\sum_{1\leq s\le a\atop s\equiv a\bmod 2}Q_s(m)
-\sum_{1\leq s\leq a\atop s\not\equiv a\bmod 2}Q_s(m+1)\,G,
\\
Q_{a-1,a}(m)=
-\sum_{1\leq s< a\atop s\equiv a\bmod 2}Q_s(m)
+\sum_{1\leq s< a\atop s\not\equiv a\bmod 2}Q_s(m+1)\,G.
\en

The above reasoning does not seem to generalize easily
to the case where more than two VO's of type I are present.

\setcounter{section}{4}
\setcounter{equation}{0}
%%%%%%%%%%%%%%%%%%%%%%%%%%%%%%%%%%%%%%%%%%%%%%%%%%%%%%%%%%%%%%%
%                                                             %
%  5. Scaling Limit                                           %
%                                                             %
%%%%%%%%%%%%%%%%%%%%%%%%%%%%%%%%%%%%%%%%%%%%%%%%%%%%%%%%%%%%%%%
\section{Scaling limit} \lb{sec:5}

Basic facts about the general RSOS models and their scaling limit 
have been worked out by Bazhanov and Reshetikhin \cite{BaRe90}, 
by diagonalizing the row-to-row transfer matrix. 
First we briefly review their results 
specializing to the present case of the ABF models in regime II. 
Then we work out new expressions for form factors 
by taking the scaling limit of the formulas on the lattice. 

%%%%%%%%%%%%%%%%%%%%%%%%%%%%
%  5.0 ==> 5.1             %
%%%%%%%%%%%%%%%%%%%%%%%%%%%%
\subsection{Review of known results}\lb{subsec:5.0}

Denote by $H$ the Hamiltonian of the lattice model
related with the row-to-row transfer matrix $T$ as 
\be
H=\frac{1}{4\pi \delta}\frac{d}{du}\log T(u)\Bigl|_{u=0} \,,
\en
where $\delta$ is the lattice spacing. 
The excitation spectrum over the ground state 
found from the Bethe ansatz has the form 
\bea
&&\epsilon_a(v)=\frac{1}{4\pi \delta}\frac{d}{dv}\log
\frac{[v+\frac{a}{2}]^*}{[v-\frac{a}{2}]^*},  
\qquad (1\le a\le k-1).
\lb{DiRe}\ena
The energies $|\epsilon_a(v)|$ 
are periodic functions with the period $-\frac{i\pi}{\log x}$, 
and have a minimum at $v=-\frac{i\pi}{2\log x}$. 
The analysis of this function shows 
that there is a non-zero gap in the spectrum, 
and the corresponding continuous 
model is a massive field theory. 

Introduce the rapidity variable $\beta$ by 
\bea
v=\frac{ik}{2\pi}\beta -\frac{i\pi}{2\log x}\,,  
\lb{Rapid}
\ena
and let 
\bea
p=e^{\frac{\pi^2}{(k+2) \log x}}
\lb{EqP}
\ena
be the temperature parameter of the lattice model. 
In the scaling limit 
\bea
&&\log x\rightarrow 0,\quad \delta\rightarrow 0,
\quad p^{\frac{k+2}{k}}\delta^{-1}\rightarrow \mbox{\it const}
=\frac{kM}{2\sin\frac{\pi}{k}}\,, 
\lb{Lim}
\ena
while keeping $\beta$ fixed, 
\eqref{DiRe} gives a massive relativistic spectrum 
\be
\lim_{\log x\rightarrow 0}
\epsilon_a\Bigl(\frac{ik}{2\pi}\beta-\frac{i\pi}{\log x}\Bigr)
=M_a\cosh\beta.
\en
Here the mass $M_a$ of the particle $a$ is defined by 
\bea
&&
M_a=\frac{\sin\frac{\pi a}{k}}{\sin\frac{\pi }{k}}\,M.
%\lb{MasRat}\\
\nn
\ena

The scattering matrices $S_{ab}(\beta)$ of these particles 
are diagonal and have a very simple form. 
For the fundamental particle/antiparticle, 
it is given by `minimal' $S$ matrices of the $\Z_k$ symmetric 
model proposed for the first time in \cite{KoSw},  
\bea
&&
S_{11}(\beta)=S_{\bar{1}\bar{1}}(\beta)
=
\frac{\sinh(\frac{\beta}{2}+\frac{i\pi}{k})}
{\sinh(\frac{\beta}{2}-\frac{i\pi}{k})},
\lb{Smin}\\
&&
S_{1\bar{1}}(\beta)=S_{\bar{1}1}(\beta)=S_{11}(i\pi -\beta). 
\nn
\ena
In general, we have 
\bea
&&S_{ab}(\beta)=
f_{{a+b}}(\beta)f_{{|a-b|}}(\beta)
\prod_{s=1}^{{\rm min}(a,b)-1} f_{{|a-b|+2s}}(\beta)^2, 
\lb{Sab}
\ena
where 
\be
&&f_{A}(\beta)={{\rm sinh}({\beta\over 2}+{i\pi A\over 2k})\over 
{\rm sinh}({\beta\over 2}-{i\pi A\over 2k})}.
\en
The particles in the scaling theory 
are not self-conjugate except for the case $k=2$ of the Ising model. 
The charge conjugation identifies 
the antiparticle $\bar{a}$ of $a$ with the particle $k-a$, 
reflecting the additional $\Z_2$ symmetry of the model.

It has been found in \cite{Za88,Ts88,BaRe90,Fat91,KlMe90} 
that the ultraviolet properties of the 
scaling theory are described by the parafermionic CFT 
\cite{ZaFa85} with the central charge 
\be
c=2{k-1\over k+2} \ .
\en
According to these works, one can treat the model 
in the scaling region as the CFT perturbed by the first energy operator 
with the left and right conformal dimensions 
\be
\Delta_{0,2}= {2\over k+2} \, .
\en

All these results agree well with 
the algebraic picture we have followed.   
The vertex operators of type II 
are identified with operators that create eigenstates of the 
row-to-row transfer matrix \cite{JM}. 
Writing the formal limit \eqref{Lim} of $\Psi_a(v)$ as $\mathcal{Z}_a(\beta)$,
where $v$ and $\beta$ are related as in \eqref{Rapid}, 
we find that the commutation relations of Proposition \ref{prop:3.2} 
and \eqref{PsiaPsib} become 
\bea 
&&\mathcal{Z}_a(\beta_1)\mathcal{Z}_b(\beta_2)
=S_{ab}(\beta_1-\beta_2)\,
\mathcal{Z}_b(\beta_2)\mathcal{Z}_a(\beta_1)\,.
\lb{ZF}
\ena
Thus the operators $\mathcal{Z}_a$ can be interpreted as generators of 
the Zamo\-lod\-chi\-kov-Faddeev algebra in the angular quantization approach 
\cite{Luk95,BrLu98}. 
The commutation relations 
between $\Phi_\pm(0)$ and $\Psi_a(v)$ (Proposition \ref{prop:3.3}) 
ensure that the states created 
by acting with the type II operators on the vacuum have the 
eigenvalues \eqref{DiRe} of the corresponding Hamiltonian.  

Our present aim is to take the continuous limit of the formulas 
in the previous section, and to interpret them as 
$\Z_k$-neutral form factors of some operators in the 
$\Z_k$-symmetric model.
We will however discuss neither the Lagrangian description of 
the present model nor the problems of identification and normalization 
of local operators. 
        
General aspects of the form factors in 
diagonal scattering theories with $\Z_k$ symmetry 
were discussed in \cite{Oota96}, where a recursive  
system of functional equations have been written. 
Another well-known example of theories with such a symmetry 
is the affine $A_{k-1}$ Toda models. 
Its $S$ matrix for fundamental particles differs from 
$S_{11}$ \eqref{Smin} by a coupling-dependent factor 
which has no poles and zeros in the physical region. 
For the affine Toda models, 
the two- and four-particle form factors have been 
determined in \cite{Oota96, Luk97b}. 
Our results below 
are very similar to those for the affine Toda case. 

%%%%%%%%%%%%%%%%%%%%%%%%%%%%
%  5.1 ==> 5.2             %
%%%%%%%%%%%%%%%%%%%%%%%%%%%%
\subsection{Projection operators}\lb{subsec:5.1}

Let us denote by $\ket{0}_m$ the $2k$ degenerate ground states in regime II.  
They are related to each other by a spatial translation. 
For the discussion of the continuous limit 
it is more convenient to deal with eigenstates 
with respect to the translation operator on the lattice,  
\bea
\frac{1}{\sqrt{2k}}
\sum_{m'=-k}^{k-1}e^{\frac{i\pi mm'}{k}}\,\ket{0}_{m'}\,.
\lb{vacm}
\ena
Let ${\rm Pr}_a$ stand for the projection operator onto the sector 
where the central height takes the value $a$. 
We shall focus attention to form factors of the 
linear transform 
\bea
\widehat{\rm Pr}_a=\frac{-1}{\sqrt{2k}}
\sum_{a'=1}^{k+1} \frac{\sin\frac{\pi a a'}{k+2}}
{\sin\frac{\pi a'}{k+2}}\,{\rm Pr}_{a'}
%\lb{TranL}
\nn
\ena
with respect to the translationally invariant vacuum $m=0$ in \eqref{vacm}.  
Thus the quantities we study in the continuous limit are  
\bea
\widehat{\mathcal{Q}}_a^{(n,n)}=
\frac{-1}{2k}
\sum_{a'=1}^{k+1}\sum_{m=-k}^{k-1}
\frac{\sin\frac{\pi a a'}{k+2}}{\sin\frac{\pi a'}{k+2}}\,
Q_{a'}^{(n,n)}(m),
\lb{fou}
\ena
where $Q_{a'}^{(n,n)}(m)$ is given by \eqref{Qbar}. 
In the angular quantization approach to affine Toda field theory, 
similar operators are 
attributed to the 
exponential operators \cite{Luk97b}. 
Note that $\widehat{\mathcal{Q}}_a^{(n,n)}=0$ if $a$ is even. 

%%%%%%%%%%%%%%%%%%%%%%%%%%%%
%  5.2 ==> 5.3             %
%%%%%%%%%%%%%%%%%%%%%%%%%%%%
\subsection{Two-particle form factors}\lb{subsec:5.2}

To illustrate 
the procedure of taking the scaling limit, 
let us consider in some details 
the simplest case $n=1$, corresponding to 
the particle-antiparticle form factors.

To find the limit \eqref{Lim} of $\widehat{\mathcal{Q}}_a^{(1,1)}$ 
it is convenient to introduce the conjugate modulus transformation.   
The standard formulas for the theta functions give 
\bea
&&
[u]=\sqrt{\frac{ik}{\tau (k+2)}}e^{-i\frac{k+2}{4k}\pi\tau}
\theta_1\Bigl(\frac{\pi u }{k+2};-\frac{k}{\tau(k+2)}\Bigr),
\nn
\\
&&
[u]^*=\sqrt{\frac{i}{\tau }}e^{-\frac{i\pi}{4}\tau}
\theta_1\Bigl(\frac{\pi u }{k};-\frac{1}{\tau}\Bigr),
\nn
\\
&&
\theta_1(u;\tau)=2\sum_{n=1}^{\infty}(-1)^{n}
e^{i\pi  \tau(n-\frac{1}{2})^2}\sin(2n-1)u.
\nn
\ena
To deal with the function $\Gamma_{m,l}$   
we borrow the technique of \cite{Jay90}, 
as mentioned in Appendix \ref{sec:app3}. 
With the notation 
$\beta=\beta_1-\beta'_1$, 
the expression \eqref{Qlat} can be rewritten in the following form.   
\begin{eqnarray}
&&
Q_a^{(1,1)}(m)=
\frac{e^{\frac{i\pi}{2k}\tau(\frac{i\beta}{\pi}+1)^2}}{k+2}
\eta(-\sfrac{1}{\tau})^2
\frac{\theta_1(\frac{\pi a}{k+2};-\frac{k}{\tau(k+2)})}
{\theta_1(\frac{\pi }{k};-\frac{1}{\tau})}
\frac{F(0)}
{F\bigl(\frac{k}{2\pi i}(\beta-i\pi)\bigr)}
\cr
&&
\quad\times
\sum_{\mu,\nu=\pm 1}\mu\nu \,
\frac{\theta_1(\frac{i\beta}{2}+\frac{\pi}{2k}(\mu+\nu)+\frac{\pi }{2};
-\frac{1}{\tau})}
{\theta_1(\frac{i\beta}{2}+\frac{\pi }{2};
-\frac{1}{\tau})}
%\sum_{l'=0}^{k+1}\sum_{m'=-k}^{k-1}
\,2\sum_{l'=0}^{k}\sum_{m'=0}^{k-1}
\sin\frac{\pi a(l'+1)}{k+2}
\cr
&&
\quad\times
e^{-\frac{i\pi mm'}{k}}
\Gamma_{m',l'}\Bigl(\frac{\mu-\nu}{2k},
-\frac{i\beta}{\pi k}-\frac{\mu+\nu}{2k}\Bigl|
-\frac{1}{\tau}\Bigr) \ .
\nn
\end{eqnarray}

{}From the sum over $m$ and $a'$, the linear transform \eqref{fou} selects 
only one term ($m'=0,l'=a-1$).  
Thus we have to analyze the limiting behavior of $\Gamma_{0,l}$, 
\bea
&&
\Gamma_{0,l}\Bigl(\frac{\mu-\nu}{2k},
-\frac{i\beta}{\pi k}-\frac{\mu+\nu}{2k}\Bigl|
-\frac{1}{\tau}
\Bigr)
\nn\\
&=\!\!&
\frac{1}{\eta(-\frac{1}{\tau})^2}
\Bigl(\sum_{n_1\geq |n_2|}-\sum_{-n_1> |n_2|}\Bigr)
(-1)^{2n_1}
p^{2\frac{k+2}{k}\bigl(
\frac{(l+1+2n_1(k+2))^2}{4(k+2)}-n_2^2k\bigr)
} 
\nn\\
&&
\qquad\qquad\times
e^{\frac{\pi i}{2k}(\mu-\nu)(l+1+2(k+2)n_1)+2\pi i  
(\frac{i\beta}{\pi }+\frac{\mu+\nu}{2})n_2}
\ ,\nn
\ena
where $p$ is defined in \eqref{EqP}
and the sum is taken over $n_1,n_2\in\frac{1}{2}\Z$ with $n_1-n_2\in\Z$. 
The leading contribution to the 
sums comes from the term with $n_1=n_2=0$. 
In the limit \eqref{Lim}, we thus find 
\be
\Gamma_{0,l}\Bigl(\frac{\mu-\nu}{2k},
-\frac{i\beta}{\pi k}-\frac{\mu+\nu}{2k}\Bigl|
-\frac{1}{\tau}
\Bigr)
= p^{2\frac{k+2}{k}(\Delta_{0,l}-\frac{c}{24})}
\Bigl(e^{\frac{\pi i }{2k}(\mu-\nu)(l+1)}
+O(p^{2\frac{k+2}{k}})\Bigr).
\nn
\en
Here $\Delta_{0,l}$ stands for 
the conformal dimension of $\Z_k$-neutral primary fields 
in the para\-fer\-mio\-nic CFT
\be
&&\Delta_{0,l}=\frac{(l+1)^2-1}{4(k+2)}
\qquad (l\equiv 0\bmod 2).
\en
 
The limit of the remaining terms can be computed directly.  
For example, 
\be
&&\lim_{\log x\rightarrow 0}
\frac{F(0)}{F\bigl(\frac{k}{2\pi i}(\beta-\pi i)\bigr)}=
F^{min}_{1\bar{1}}(\beta)\,,
\en 
where we introduced the minimal form factor for the
particle-antiparticle scattering \cite{KaWe78} 
\be
&&F_{1\bar{1}}^{min}(\beta)
=\prod_{n=1}^{\infty}
\frac{\Gamma(\frac{i\beta}{2\pi}+\frac{1}{k}+n+\frac{1}{2})
\Gamma(-\frac{i\beta}{2\pi}+\frac{1}{k}+n-\frac{1}{2})
\Gamma(n-\frac{1}{k})^2}
{\Gamma(\frac{i\beta}{2\pi}-\frac{1}{k}+n+\frac{1}{2})
\Gamma(-\frac{i\beta}{2\pi}-\frac{1}{k}+n-\frac{1}{2})
\Gamma(n+\frac{1}{k})^2}
\ ,
\en
normalized as $F^{min}_{1\bar{1}}(i\pi)=1$. 
As for the sum over different components of vertex operators we use
\be
&&
\sum_{\mu,\nu=\pm 1}\mu\nu\, e^{\frac{\pi i }{2k}(\mu-\nu)a}             
\cosh\Bigl(\frac{\beta}{2}+\frac{i\pi}{2k}(\mu+\nu)\Bigr)
\\
&&
=
4\cosh\frac{\beta}{2}\;\sin\frac{\pi(a-1)}{2k}\;\sin\frac{\pi(a+1)}{2k}.
\en
This leads us to 
the simple expression for the two-point form factor 
of the projection operators 
\be
  &&\widehat{\mathcal{Q}}_a^{(1,1)}
  =-2\,\frac{\sin\frac{\pi(a-1)}{2k}\sin\frac{\pi(a+1)}{2k}}{\sin\frac{\pi}{k}}
  \,p^{2\frac{k+2}{k}\Delta_{0,a-1}}F^{min}_{1\bar{1}}(\beta) 
  \times\bigl(1+o(1)\bigr).
\en
As we expect from the results of \cite{Oota96}, 
there is no singularity at $\beta=i\pi$.  

In the next subsection, we compute the scaling limit 
for the case of multi-particles in a similar manner. 
We note in passing that 
\be
&&\lim_{\log x\rightarrow 0}
\frac{F\bigl(\frac{k}{2\pi i}\beta\bigr)}{F(0)}=
\frac{F^{min}_{11}(\beta)}
{\sinh\frac{\beta}{2}\sinh(\frac{\beta}{2}+\frac{i\pi}{k})
}\ ,
\en
where 
the minimal two particle form factor $F^{min}_{11}(\beta)$ reads \cite{KaWe78}
\be
&&F^{min}_{11}(\beta)=
\sinh\frac{\beta}{2}\sinh\Bigl(\frac{\beta}{2}+\frac{i\pi}{k}\Bigr)
\\
&&\qquad \quad\times 
\prod_{n=1}^{\infty}
\frac{\Gamma(\frac{i\beta}{2\pi}-\frac{1}{k}+n)
\Gamma(-\frac{i\beta}{2\pi}-\frac{1}{k}+n)\Gamma(n+\frac{1}{k})^2}
{\Gamma(\frac{i\beta}{2\pi}+\frac{1}{k}+n)
\Gamma(-\frac{i\beta}{2\pi}+\frac{1}{k}+n)\Gamma(n-\frac{1}{k})^2}.
\en

%%%%%%%%%%%%%%%%%%%%%%%%%%%%
%  5.3 ==> 5.4             %
%%%%%%%%%%%%%%%%%%%%%%%%%%%%
\subsection{Many particles}\lb{subsec:5.3}

To present the result in the general case, it is convenient to 
use the formal bosonization rule as explained in \cite{Luk97b,Luk97a}.  
Let us introduce 
\bea
&&\mathcal{B}^a_{1}(\beta)
=-\sum_{\mu=\pm 1}\mu\,e^{-\frac{\mu \pi ia}{2k}}\mathcal{Z}_{1,\mu}(\beta),
\lb{expo1}\\
&&\mathcal{B}^a_{\bar{1}}(\beta)
=\sum_{\mu=\pm 1}\mu\,e^{\frac{\mu \pi ia}{2k}}
\mathcal{Z}_{\bar{1},\mu}(\beta).
\lb{expo2}
\ena
Here $\mathcal{Z}_{1,\mu}(\beta)$ and $\mathcal{Z}_{\bar{1},\mu}(\beta)$ 
are some operators,    
for which we assume the contraction rules 
\be
&&
\dbr{\mathcal{Z}_{1,\mu}(\beta_1)\mathcal{Z}_{1,\nu}(\beta_2)}
=\dbr{\mathcal{Z}_{\bar{1},-\mu}(\beta_1)\mathcal{Z}_{\bar{1},-\nu}(\beta_2)}
\\
&&
\qquad=
\frac{F^{min}_{11}(\beta)}
{\sinh(\frac{\beta}{2}+\frac{i\pi}{k})\sinh(\frac{\beta}{2}-\frac{i\pi}{k})}
\frac{\sinh\bigl(\frac{\beta}{2}-\frac{i\pi}{2k}(\mu-\nu)\bigr)}
{2\sin\frac{\pi}{k}\:\sinh\frac{\beta}{2}},
\\
&&
\dbr{\mathcal{Z}_{1,\mu}(\beta_1)\mathcal{Z}_{\bar{1},\nu}(\beta_2)}
=\dbr{\mathcal{Z}_{\bar{1},-\mu}(\beta_1)\mathcal{Z}_{1,-\nu}(\beta_2)}
\\
&&
\qquad =
F^{min}_{1\bar{1}}(\beta)
\frac{\cosh\bigl(\frac{\beta}{2}+\frac{i\pi}{2k}(\mu+\nu)\bigr)}
{2\sin\frac{\pi}{k}\:\cosh\frac{\beta}{2}},
\en
where $\beta=\beta_1-\beta_2$, 
and the Wick theorem as in \eqref{Wick} where $C$
is replaced by 
\be
C_1=\frac{1}{\sqrt{2\sin\frac{\pi}{k}}}\,.
\en
The prescription \eqref{expo1}, \eqref{expo2}
is analogous to the one for 
the free field representation of 
form factors of the `exponential fields' \cite{Luk97a, Luk97b} 
in the sine-Gordon and affine Toda theories. 
The exponential factors $e^{\pm\frac{\pi ia}{2k}}$   
are the only remnant of the `zero mode' part. 

The limit of the formula
\eqref{fou} can be written compactly as 
\bea
  &&\quad
  \lim_{\log x \rightarrow 0}\Bigl(
  p^{-2\frac{k+2}{k}\Delta_{0,a-1}}\widehat{\mathcal{Q}}_a^{(n,n)}\Bigr)
  \lb{answer}\\
  &&=\dbr{\mathcal{B}^{a}_{\bar{1}}(\beta_1)\cdots
  \mathcal{B}^{a}_{\bar{1}}(\beta_n)
  \mathcal{B}^{a}_{1}(\beta'_1)\cdots\mathcal{B}^{a}_{1}(\beta'_n)}
  \nn\\
  &&=
  \prod_{1\le i<j\le n}\!\!\frac{F^{min}_{11}(\beta_i-\beta_j)}
         {(x_i-\omega^2 x_j)(x_i-\omega^{-2} x_j)}\cdot\!\!
  \prod_{1\le i<j\le n}\!\!\frac{F^{min}_{11}(\beta'_i-\beta'_j)}
         {(y_i-\omega^2 y_j)(y_i-\omega^{-2} y_j)}
  \nn\\
  &&\quad\times
  \prod_{1\le i, j\le n}\!\!\frac{F^{min}_{1\bar{1}}(\beta_i-\beta'_j)}
         {x_i+y_j}
  \times C_1^{2n}\,2^{2n(n-1)}\sigma_n^{n-1}\tau_n^{n-1}
  R_{\frac{a+1}{2}}^{(n,n)}(x\,;y).
  \nn
\ena
The notation is as follows. 
We set $x_i=e^{\beta_i}$, $y_i=e^{\beta'_i}$, 
$\sigma_r=\sigma_r(x)$, $\tau_r=\sigma_r(y)$ with 
\be
\prod_{j=1}^n(t+x_j)=\sum_{r=0}^n t^{n-r}\sigma_r(x),
\en
$\omega=e^{\pi i/k}$, 
and $\{l\}=\omega^l-\omega^{-l}$. 
The polynomials $R_\alpha^{(m,n)}(x\,;y)$ are defined by
\be
&&
R_\alpha^{(m,n)}(x\,;y)
=\sum_{\mu_1,\cdots,\mu_m=\pm 1 \atop \nu_1,\cdots,\nu_n=\pm 1 }
\prod_{i=1}^m\mu_i\omega^{(\alpha-m+n)\mu_i}\cdot
\prod_{j=1}^n \nu_j\omega^{(\alpha+m-n)\nu_j}
\\
&&
\quad\times
\prod_{1\le i<j\le m}\!\!
\frac{x_i\omega^{\mu_i}-x_j\omega^{\mu_j}}{x_i-x_j}\cdot\!\!
\prod_{1\le i<j\le n}\!\!
\frac{y_i\omega^{\nu_i}-y_j\omega^{\nu_j}}{y_i-y_j}\cdot\!\!
\prod_{1\le i\le m\atop 1\le j\le n}\!\!
(x_i\omega^{-\mu_i}+y_j\omega^{-\nu_j}).
\en
For example, 
\be
&&
R_\alpha^{(1,1)}=\{\alpha\}\{\alpha-1\}(\sigma_1+\tau_1),
\\&&
R_\alpha^{(2,2)}=\{\alpha\}\{\alpha-1\}
\bigl(\{\alpha\}\{\alpha-1\}
(\sigma_1+\tau_1)(\sigma_2\tau_1+\sigma_1\tau_2)
\\
&&\qquad\qquad\qquad\qquad\quad
+\{\alpha+1\}\{\alpha-2\}(\sigma_2-\tau_2)^2\bigr),
\\
&&R^{(3,3)}_\alpha=\{\alpha\}^3\{\alpha-1\}^3
(\sigma_1+\tau_1)
(\sigma_3+\sigma_2\tau_1+ \sigma_1\tau_2+\tau_3)
(\sigma_3\tau_2+ \sigma_2\tau_3)\\
&&\quad-
 \{\alpha+2\}\{\alpha\}^2\{\alpha-1\}^2\{\alpha-3\}
(\sigma_1+\tau_1)(\sigma_3+\tau_3)(\sigma_3\tau_2+ \sigma_2\tau_3)\\
&&\quad+
  \{\alpha+1\}^2\{\alpha\}\{\alpha-1\}\{\alpha-2\}^2
(\sigma_2\tau_1+\sigma_1\tau_2)
(\sigma_3+\tau_3)^2\\
&&\quad+
 \{\alpha+1\}\{\alpha\}^2\{\alpha-1\}^2\{\alpha-2\}
 \bigl(
 (\sigma_1+\tau_1) (\sigma_3\tau_1-\sigma_1\tau_3)^2\\
 &&\qquad\quad
 +
 (\sigma_2-\tau_2)^2 (\sigma_3\tau_2+\sigma_2\tau_3)
 -3 
 (\sigma_1+\tau_1)(\sigma_3+\tau_3)(\sigma_3\tau_2+ \sigma_2\tau_3) \bigr)\\
&&\quad+
\{\alpha+2\} \{\alpha+1\}\{\alpha\}\{\alpha-1\}\{\alpha-2\}\{\alpha-3\}
(\sigma_3+\tau_3)^3.
\en
We note that $R_{1}^{(n,n)}=0$ holds for all $n$. 

It can be shown that $R^{(m,n)}_\alpha(x\,;y)$ is a sum of products of  
determinants.
Namely let $\Lambda(m,n)$ denote the set of partitions 
$\lambda=(\lambda_1,\cdots,\lambda_m)$ 
satisfying $n\ge \lambda_1\ge\cdots\ge\lambda_m\ge 0$. 
For $\lambda\in\Lambda(m,n)$ we write 
$\widetilde{\lambda'}=(m-\lambda'_n,\cdots,m-\lambda'_1)$, 
where $\lambda'=(\lambda'_1,\cdots,\lambda'_n)$ denotes the conjugate
partition. 
Then 
\bea
R^{(m,n)}_\alpha(x\,;y)&\!\!=\!\!&
\prod_{i=1}^m\{\alpha-i+n\}\cdot
\prod_{j=1}^n\{\alpha-j+m\}
%\lb{Schur}\\
\nn\\
&&\times 
\sum_{\lambda \in \Lambda(m,n)}
S_\lambda(x\,;\alpha+n,\omega) S_{\widetilde{\lambda'}}(y\,;\alpha+m,\omega),
\nn
\ena
where
\begin{eqnarray*}
&&S_\lambda(x\,;\alpha,\omega)=
\det \left(
\frac{\{\alpha-\lambda'_i+i-2j\}}{\{\alpha-i\}}
\sigma_{\lambda'_i-i+j}(x)\right)_{1\leq i,j\leq N},  
\end{eqnarray*}
with $N\ge n$.

\setcounter{section}{5}
\setcounter{equation}{0}
%%%%%%%%%%%%%%%%%%%%%%%%%%%%%%%%%%%%%%%%%%%%%%%%%%%%%%%%%%%%%%%
%                                                             %
%  6. Discussions                                             %
%                                                             %
%%%%%%%%%%%%%%%%%%%%%%%%%%%%%%%%%%%%%%%%%%%%%%%%%%%%%%%%%%%%%%%
\section{Discussions}

In this paper we have developed
an algebraic approach to the ABF models in regime II, 
using free fields. 
We have found the following. 

First, our results immediately give an integral representation for 
correlation functions of these non-critical models. 
We computed exactly the simplest integrals for the nearest neighbor 
correlation functions. 
The result supports the validity of our construction. 
Beyond this case, it remains a technical open problem to perform
multiple contour integrals explicitly.

Second, following the prescription of \cite{JM} we found states
which diagonalize the row-to-row transfer matrix
in the thermodynamic limit.
These are obtained by acting with type II vertex operators on the vacuum. 
Compared with other lattice models studied so far, 
a nice feature is that these operators do not contain contour integrals. 
This gives us a hope for handling them more effectively in  
further analysis.  

As an application we studied the continuous limit of 
the traces of vertex operators of type II. 
We obtained a family of functions satisfying 
Watson's equations as well as kinematical pole conditions. 
The former is a consequence of the commutation relation
\eqref{ZF} and the property of  the corner Hamiltonian as the grading operator,
while the latter follows
from
\be
&&
i\mathcal{B}^a_{\bar{1}}(\beta_1)\mathcal{B}^a_{1}(\beta_2)
=
\frac{1}{\beta_1-\beta_2-i\pi}\times \id +O(1)
\quad (\beta_1\rightarrow \beta_2+i\pi).
\en
On these grounds, 
we interpreted the resulting formulas as 
form-factors of some local operators 
in the $\Z_k$ invariant massive scattering theory with minimal $S$-matrices. 
Our analysis is not complete, since we have considered only 
the $\Z_k$ neutral sector, which corresponds to having 
the same number of $\Psi_+(v)$'s and $\Psi_-(v)$'s in the trace.  
The bound state conditions also remain to be worked out. 

There are other unsolved problems which deserve attention.
Of particular interest is 
the problem of identifying the local operators corresponding to
the projection operators. 
The form-factors $\widehat{\mathcal{Q}}_a^{(n,n)}$ 
are proportional to $M^{2\Delta_{0,a-1}}$. 
This mass dependence indicates that the projection operators
correspond to perturbations of the neutral primary fields 
in the parafermionic CFT with the respective conformal 
dimensions.\footnote{A similar study 
has been made in \cite{La97} for the scaling limit of the XXZ model.}
An argument in favor of this proposal is that 
the cluster property characteristic to the `exponential operators' 
seems to hold. 

It is natural to identify the case 
$a=1$, $\Delta_{0,0}=0$, with the form-factor of the identity operator.
The fact that $\widehat{\mathcal{Q}}_1^{(n,n)}=0$ for all $n$
agrees with this identification.  
In the case  $a=3$, $\Delta_{0,2}=\frac{2}{k+2}$, 
we infer that the $\widehat{\mathcal{Q}}_3^{(n,n)}$ give the 
form factors of the first energy operator. 
Since this operator perturbs the parafermionic CFT into the massive region, 
the corresponding form factors are
proportional to those of the trace of stress-energy tensor. 
One can show that, for $a=3$, the polynomial $R^{(n,n)}_{2}(x;y)$ in the 
formula \eqref{answer} contains a factor 
\be
(\sigma_1+\tau_1)
\Bigl(\frac{\sigma_{n-1}}{\sigma_n}+\frac{\tau_{n-1}}{\tau_{n}}\Bigr), 
\en
as we expect for form factors of the stress-energy tensors 
\cite{Zam91,Oota96}. 

However, our understanding of this 
problem of identification of local operators, 
as well as that of deriving the 
vacuum expectation values of local fields, is still incomplete.

In the scaling limit we 
considered only the translationally invariant sector $m=0$. 
{}From the lattice point of view, the other sectors have equal rights. 
It would be interesting to understand the field theoretical meaning of 
the traces corresponding to the $m\neq 0$ sectors. 

Although practically all the technical elements of the construction 
have been known, 
we think that the interpretation as a free field realization 
of ABF models in regime II is rather interesting. 
It opens up a way for applying algebraic methods to the study of
physical phenomena in these models.

\bigskip

%%%%%%%%%%%%%%%%%%%%%%%%%%%%%%%%%%%%%%%%%%%%%%%%%%%%%%%%%%%%%%%
%                                                             %
%  Acknowledgments                                            %
%                                                             %
%%%%%%%%%%%%%%%%%%%%%%%%%%%%%%%%%%%%%%%%%%%%%%%%%%%%%%%%%%%%%%%
\noindent
{\it Acknowledgments.}\quad 

We thank 
Michael Lashkevich and Takeshi Oota for helpful discussions. 
M.~J. is grateful to Professor Shi-Shyr Roan  
for kind invitation and warm hospitality during his visit to 
Academia Sinica, Taiwan. 
H.~K. thanks V.~A.~Fateev and Al.~B.~Zamolodchikov for discussions and 
colleagues in LAPTH for kind hospitality during his stay in Annecy. 
This work is partially supported by  
the Grant-in-Aid for Scientific Research 
no.10440042 and no.11640030, the Ministry of Education. 

\appendix
\setcounter{equation}{0}
%%%%%%%%%%%%%%%%%%%%%%%%%%%%%%%%%%%%%%%%%%%%%%%%%%%%%%%%%%%%%%%
%                                                             %
%  App.A. Formulas for operators                              %
%                                                             %
%%%%%%%%%%%%%%%%%%%%%%%%%%%%%%%%%%%%%%%%%%%%%%%%%%%%%%%%%%%%%%%
\section{Formulas for operators}\lb{sec:app1}

We give here explicit formulas for the operators in  
the deformed parafermion theory. 
\medskip

\noindent {\bf Oscillators}\quad
We consider two kinds of oscillators $a_{j,n}$ 
($j=1,2$, $n\in\Z\backslash\{0\}$) 
and `zero-mode' operators $P_j,Q_j$ ($j=1,2$) satisfying 
\begin{eqnarray*}
&& [a_{1,n},a_{1,n'}]=\frac{[2n]_x[(k+2)n]_x}{n}\delta_{n+n',0},
\\
&& [a_{2,n},a_{2,n'}]=-\frac{[2n]_x[kn]_x}{n}\delta_{n+n',0},
\\
&&[P_1,Q_1]=2(k+2),\quad [P_2,Q_2]=-2k, 
\end{eqnarray*}
where 
\be
[n]_x=\frac{x^n-x^{-n}}{x-x^{-1}}.
\en

\noindent {\bf Fock space}\quad
For $l,m\in\Z$ we set
\be
&&\F_{m,l}=\C[a_{1,-n},a_{2,-n}~~(n>0)]\,\ket{m,l},
\\
&&
P_1\ket{m,l}=l\ket{m,l},
\quad
P_2\ket{m,l}=-m\ket{m,l},
\\
&&e^{\frac{Q_1}{2(k+2)}}\ket{m,l}=\ket{m,l+1},
\quad
e^{\frac{Q_2}{2k}}\ket{m,l}=\ket{m+1,l}.
\end{eqnarray*}

\noindent {\bf Grading operator}\quad
We set $D=D^{osc}+D^{zero}$, with 
\begin{eqnarray*}
&& D^{osc}=\sum_{n>0}\frac{n^2}{[2n]_x[(k+2)n]_x}a_{1,-n}a_{1,n}
-\sum_{n>0}\frac{n^2}{[2n]_x[kn]_x}a_{2,-n}a_{2,n},
\\
&&D^{zero}=\frac{(P_1+1)^2-1}{4(k+2)}-\frac{P_2^2}{4k}.
\end{eqnarray*}
We have 
\be
&&[D,a_{j,n}]=-n a_{j,n}, \\
&&[D,Q_1]=P_1+1, \quad [D,Q_2]=P_2, \\
&&
D\ket{m,l}
=\Bigl(\frac{(l+1)^2-1}{4(k+2)}-\frac{m^2}{4k}\Bigr)\ket{m,l}. 
\en

\noindent {\bf Notational convention}\quad 
We use the following notation. 
\be
&&\phi_j(A;B,C|z;C')=-\frac{A}{BC}\left(Q_j+P_j\log z\right)
+\sum_{n\neq 0}\frac{[An]_x}{[Bn]_x[Cn]_x}a_{j,n}z^{-n}x^{C'|n|},
\\
&&\phi_j(B|z;C)=\phi_j(A;A,B|z;C),
\\
&&
\phi_j^{(\pm)}(A;B|z)=\frac{A}{B}P_j\log x
+(x-x^{-1})\sum_{n> 0}\frac{[An]_x}{[Bn]_x}a_{j,\pm n}z^{\mp n}.
\en
We use the `additive' parameters $u,u',\cdots$ and  
the `multiplicative' ones $z,z',\cdots$ on an equal footing.  
Unless otherwise stated explicitly, they are related by 
\[
z=x^{2u},\quad z'=x^{2u'},\cdots.
\]

The following are the list of operators used in the text. 
\medskip

\noindent {\bf Type II VO's}

\begin{eqnarray*}
&&\Psi_\pm(u)=\mp\frac{1}{x-x^{-1}}\left(
\Psi_{\pm,+}(u)-\Psi_{\pm,-}(u)\right),
\\
&&\Psi_{\pm,\ve}(u)=z^{-\frac{1}{k}}\times
\\
&&\quad :\exp\Bigl(
\pm \phi_2(k|z;\pm \frac{k}{2})
\pm\ve\phi_1^{(\ve)}(1;2|zx^{\pm\varepsilon\frac{k}{2}})
-\ve\phi_2^{(\ve)}(1;2|zx^{\pm\varepsilon \frac{k+2}{2}})
\Bigr):.
\end{eqnarray*}
More generally we set
\bea
&&\Psi_a(u)=z^{-\frac{a(a-1)}{k}}\sum_{j=0}^a(-1)^j \,x^{(a-1)(j-\frac{a}{2})}
\prod_{i=1}^j\frac{[a-i+1]_x}{[i]_x}
\lb{fusP}\\
&&\quad\times
:\prod_{i=1}^{a-j}\Psi_{-,-}(-\sfrac{a+1}{2}+i+u)\cdot\!\!\!
\prod_{i=a-j+1}^{a}\Psi_{-,+}(-\sfrac{a+1}{2}+i+u):.
\nn
\ena

\noindent {\bf Screening current}\quad 
\begin{eqnarray*}
&&S(u)=-\frac{1}{x -x^{-1}}\left(S_+(u)-S_-(u)\right),
\\
&&
S_\pm(u)=
\\
&&\quad:\exp\Bigl(\phi_1(k+2|z;-\frac{k+2}{2})
\pm\phi^{(\pm)}_2(1;2|zx^{\mp\frac{k+2}{2}})
\pm\phi^{(\pm)}_1(1;2|zx^{\mp\frac{k}{2}})
\Bigr):.
\en

\noindent {\bf Type I VO's}

\be
&&\Phi_{-}(u)=z^{\frac{k-1}{2k(k+2)}}\times
\\
&&\quad 
:\exp\Bigl(
-\phi_1(1;2,k+2|zx^{k};\frac{k+2}{2})
-\phi_2(1;2,k|zx^{k};\frac{k}{2})
\Bigr):,
\\
&&\Phi_{+}(u)=
\oint_{C_\Phi(z)}\underline{ dz}' \Phi_{-}(u)
S(u')\frac{[u-u'-\frac{3}{2}-P_1]}{[u-u'-\frac{3}{2}]},
\\
&& 
\Phi^*_{-}(u)= \frac{g}{[P_1+1]}\,z^{\frac{k-1}{2k(k+2)}}\times 
\\
&&\quad 
:\exp\Bigl(
-\phi_1(1;2,k+2|z;-\frac{k+2}{2})
-\phi_2(-1;2,k|z;-\frac{k}{2})
\Bigr):,
\\
&& \Phi^*_{+}(u)=- \oint_{C_{\Phi^*}(z)}\underline{ dz}' \Phi^*_{-}(u)
S(u')\frac{[u-u'-\frac{1}{2}-P_1]}{[u-u'-\frac{1}{2}]}. 
\en
The contours $C_{\Phi}(z),C_{\Phi^*}(z)$ are specified by the rules 
\begin{eqnarray*}
C_{\Phi}(z):\quad z'&\!\!=\!\!&zx^{-3+ 2(k+2) n} \qquad(n> 0)\qquad\qquad 
     {\rm inside},\\
            &\!\!=\!\!&zx^{-1+ 2(k+2) n} \qquad(n\leq 0)\qquad\qquad 
     {\rm outside},\\
C_{\Phi^*}(z):\quad z'&\!\!=\!\!&zx^{-1+ 2(k+2) n} \qquad(n\geq 0)\qquad\qquad 
     {\rm inside},\\
            &\!\!=\!\!&zx^{1+2(k+2) n} \,\,\,\qquad(n\leq 0)\qquad\qquad 
     {\rm outside}, 
\end{eqnarray*}
and
\be
g&\!\!=\!\!&(x-x^{-1})x^{-\frac{1}{k+2}}\\
&&\times
\frac{(x^{2k+2},x^{4k+2};x^{2k},x^{2k+4})_\infty}
{(x^{2k},x^{4k+4};x^{2k},x^{2k+4})_\infty}
(x^{2k+2};x^{2k+4})_\infty^2 (x^{2k+4};x^{2k+4})_\infty.
\en

\noindent {\bf $\xi$-$\eta$ system}\quad 

\be
&&\xi(u)=\,:\exp\Bigl(
\phi_1(2|z;\frac{k}{2})+\phi_2(2|z;\frac{k+2}{2})
\Bigr):,\\
&&
\eta(u)=\,:\exp\Bigl(
-\phi_1(2|z;\frac{k}{2})-\phi_2(2|z;\frac{k+2}{2})
\Bigr):.
\end{eqnarray*}

\medskip

\noindent{\bf Deformed $W$ currents}\quad 
The first DWA current $W^1(u)$ is 
\be
  &&W^1(u)=[k+1]_x\Lambda_-(u)+[k+2]_x\Lambda_0(u)+[k+1]_x\Lambda_+(u),
\en
where $\Lambda_{\pm}(u)$ and $\Lambda_0(u)$ are
\be
  \Lambda_{\pm}(u)&\!\!=\!\!&
  x^{\pm1}g^{-1}x^{\frac{1-k}{k}}z^{\frac{1-k}{k(k+2)}}
  :\Phi_-(k+2+u)[P]\Phi^*_-(u)S_-(\sfrac{k+1}{2}\mp\sfrac{k+2}{2}+u):\,,\\
  \Lambda_0(u)&\!\!=\!\!&-g^{-1}x^{\frac{1-k}{k}}z^{\frac{1-k}{k(k+2)}}
  :\Phi_-(k+2+u)[P]\Phi^*_-(u)S_+(k+\sfrac32+u):.
\en
In general the DWA currents $W^j(u)$ ($j=1,2,\cdots$) are given by
\be
  W^j(u)&\!\!=\!\!&\!\!\sum_{a,b,c\geq 0\atop a+b+c=j}C^j_{a,b,c}
 :\prod_{i=1}^a\Lambda_-(\sfrac{j+1}{2}-i+u)\\
  &&\qquad\quad\times
  \prod_{i=a+1}^{a+b}\Lambda_0(\sfrac{j+1}{2}-i+u)\cdot\!\!\!  
  \prod_{i=a+b+1}^j\Lambda_+(\sfrac{j+1}{2}-i+u):\,,  
\en
where $C^j_{a,b,c}$ is
\be
  &&C^j_{a,b,c}=\prod_{i=1}^a\frac{[k+2-i\,]_x}{[i\,]_x}\cdot
  \prod_{i=1}^b\frac{[k+1+i\,]_x}{[i\,]_x}\cdot
  \prod_{i=1}^c\frac{[k+2-i\,]_x}{[i\,]_x}.
\en
\medskip

\noindent {\bf Comparison of notation}\quad
We have followed the notation in \cite{Konno} with minor changes.  
Referring to the symbols in \cite{Konno} by the subscript `$K$', 
we have
\be  
&&x=q_K,  \\
&&\F_{m,l}=\F^{PF}_{l,m,K}, \\
&&S(u)=S_K(z), \\
&&\eta(u)=\eta_K(z),\\
&&\Phi_{-}(u)=\phi_{1,1,K}(q^kz)z^{\frac{k-1}{2k(k+2)}}, \\
&&\Phi^*_{-}(u)=\frac{g}{[P_1+1]}
\phi_{1,-1,K}(z)z^{\frac{k-1}{2k(k+2)}}, \\
&&\Psi_+(u)=\Psi^\dagger_K(z)z^{-\frac{1}{k}}, \\
&&\Psi_-(u)=\Psi_K(z)z^{-\frac{1}{k}}.
\en
We have introduced fractional powers of $z$ in order that 
the homogeneity property 
\be
x^{2vD}Y(u)x^{-2vD}=Y(u+v)
\en
holds for $Y=S,\Phi_{\pm},\Phi_{\pm}^*,\Psi_{\pm},\xi,\Psi_a,W^j$. 
As for $\eta$ we have $x^{2vD}\eta(u)x^{-2vD}=\eta(u+v) x^{2v}$.

\setcounter{equation}{0}
%%%%%%%%%%%%%%%%%%%%%%%%%%%%%%%%%%%%%%%%%%%%%%%%%%%%%%%%%%%%%%%
%                                                             %
%  App.B. Resolution by $\xi$-$\eta$ system                   %                
%                                                             %
%%%%%%%%%%%%%%%%%%%%%%%%%%%%%%%%%%%%%%%%%%%%%%%%%%%%%%%%%%%%%%%
\section{Resolution by $\xi$-$\eta$ system}\lb{sec:app2}

We summarize the main points concerning the resolution by the 
$\xi$-$\eta$ system.

On the space $\F_{m,l}$ with $m\equiv l\bmod 2$, we have 
expansions of the form 
\be
&&\eta(u)=\sum_{n\in\Z}\eta_n z^{-n-1},\qquad
\xi(u)=\sum_{n\in\Z}\xi_n z^{-n}.
\en
The Fourier components satisfy
\begin{eqnarray*}
&&[\eta_n,\eta_{n'}]_+=0,
\quad [\xi_n,\xi_{n'}]_+=0,\quad 
[\eta_n,\xi_{n'}]_+=\delta_{n+n',0}.
\end{eqnarray*}
The components $\eta_0,\xi_0$ commute with $D$. 
Since $\xi_0^2=\eta_0^2=0$ and $\xi_0\eta_0+\eta_0\xi_0=\id$, 
the complex 
\be
&&
\cdots \overset{\eta_0}{\longrightarrow }
\F_{m-k,l-(k+2)}\overset{\eta_0}{\longrightarrow }
\F_{m,l}\overset{\eta_0}{\longrightarrow }
\F_{m+k,l+k+2} \overset{\eta_0}{\longrightarrow }
\F_{m+2k,l+2(k+2)}\overset{\eta_0}{\longrightarrow }
\cdots
\en
is exact. 
Let 
\be
\Ft_{m,l}&=&{\rm Ker}\,\eta_0\bigl|_{\F_{m,l}} 
\\
\quad &=&{\rm Coker}\,\eta_0\bigl|_{\F_{m-k,l-(k+2)}}\,.
\en
Then, for an operator $\mathcal{O}$ on $\oplus_{m\equiv l\bmod 2}\F_{m,l}$ 
commuting with $\eta_0$, we have 
\be
\tr_{\Ft_{m,l}}(\mathcal{O})
&=&
\sum_{n\ge 0}(-1)^n\,\tr_{\F_{m+kn,l+(k+2)n}}(\mathcal{O})
\\
&=&
-\sum_{n<0}(-1)^n\,\tr_{\F_{m+kn,l+(k+2)n}}(\mathcal{O})\,.
\en

\setcounter{equation}{0}
%%%%%%%%%%%%%%%%%%%%%%%%%%%%%%%%%%%%%%%%%%%%%%%%%%%%%%%%%%%%%%%
%                                                             %
%  App.C. Calculation of the trace                            %
%                                                             %
%%%%%%%%%%%%%%%%%%%%%%%%%%%%%%%%%%%%%%%%%%%%%%%%%%%%%%%%%%%%%%%
\section{Calculation of the trace}\lb{sec:app3}

We outline the computation of the trace of a product of type II
operators focusing attention to the neutral case 
\bea
&&\mathcal{O}=
x^{2kD}
\Psi_{+}(v_1)\cdots \Psi_{+}(v_n)
\Psi_{-}(v'_1)\cdots\Psi_{-}(v'_n).
\lb{ou}
\ena
It is convenient in what follows to consider the oscillator part 
and the zero mode part separately. 
Let us write 
\be
&&\Psi_{\pm,\varepsilon}(v)=\Psi^{osc}_{\pm,\varepsilon}(v)
\Psi^{zero}_{\pm,\varepsilon}(v),
\\
&&\Psi^{zero}_{\pm,\varepsilon}(v)=
e^{\mp Q_2/k}z^{\mp P_2/k-1/k}x^{\pm\varepsilon P_1/2-\varepsilon P_2/2}. 
\en
\medskip

\noindent{\bf {Oscillators}}\quad
The contributions from the oscillator part
\be
\dbr{\mathcal{O}}=
\frac{\tr_{\F_{m,l}}\bigl(x^{2kD^{osc}}\mathcal{O}^{osc}\bigr)}
{\tr_{\F_{m,l}}(x^{2kD^{osc}})}
\en
are given by the following rules: 
\bea
&&
C^{-N}\dbr{\Psi_{\varepsilon_1,\varepsilon'_1}(v_1)\cdots
\Psi_{\varepsilon_N,\varepsilon'_N}(v_N)}
\lb{Wick}\\
&&\qquad\qquad=\prod_{1\le i<j\le N}C^{-2}
\dbr{\Psi_{\varepsilon_i,\varepsilon'_i}(v_i)
\Psi_{\varepsilon_j,\varepsilon'_j}(v_j)}
\nn
\ena
for $\sum_{i=1}^N\varepsilon_i=0$, and 
\bea
  &&\;\;\dbr{\Psi_{+,\ve_1}(v_1)\Psi_{+,\ve_2}(v_2)}=
  \dbr{\Psi_{-,-\ve_1}(v_1)\Psi_{-,-\ve_2}(v_2)}
\lb{contr1}\\
  &&\qquad=C^2 F(v)\frac{[v+\frac{\ve_1-\ve_2}{2}]^*}{[v-1]^*}\,
  x^{-\frac{2}{k}(1+\frac{\ve_1-\ve_2}{2})v
  +\frac{1+\ve_1\ve_2}{2k}+1+\frac{\ve_1-\ve_2}{2}},
\nn\\
  &&\;\;\dbr{\Psi_{+,\ve_1}(v_1)\Psi_{-,\ve_2}(v_2)}=
  \dbr{\Psi_{-,-\ve_1}(v_1)\Psi_{+,-\ve_2}(v_2)}
\lb{contr2}\\
  &&\qquad=C^2 F(v-\sfrac{k}{2})^{-1}
  \frac{[v-\frac{k}{2}-\frac{\ve_1+\ve_2}{2}]^*}{[v-\frac{k}{2}]^*}\,
  x^{\frac{\ve_1+\ve_2}{k}v
  -\frac{1+k}{2k}(1+\ve_1\ve_2)-\frac{\ve_1+\ve_2}{2}},
\nn
\ena
where $v=v_2-v_1$, $F(v)$ is given in \eqref{F(w)},
and
\be
  C=(x^{2k};x^{2k})_\infty
  \frac{(x^{2+4k};x^{2k},x^{2k})_\infty}{(x^{-2+2k};x^{2k},x^{2k})_\infty}.
\en
\medskip

\noindent{\bf Zero mode}\quad 
The trace over the zero mode in the Fock space 
is taken in a standard fashion, following \eqref{resol1},\eqref{resol2}. 
As we show at the end of this appendix, 
our operator \eqref{ou} satisfies in addition 
\bea 
\tr_{\F_{m,l}}(\mathcal{O}_0)=\tr_{\F_{m,-l-2}}(\mathcal{O}_1).
\lb{physical}
\ena
Then, for any fixed real numbers $h,h'\in \R$, 
the trace can be rewritten as follows. 
\bea
&&\quad\tr_{H^0(\Cc_{m,l})}(\mathcal{O})
%\lb{Ind}\\
\nn\\
&=\!\!&\Bigl(\sum_{s\le h}\sum_{n\ge 0}-\sum_{s> h}\sum_{n< 0}\Bigr)
(-1)^n\,\tr_{\F_{m+kn,l+(k+2)(n-2s)}}(\mathcal{O}_0)
\nn\\
&&
-\Bigl(\sum_{s<-h'}\sum_{n\ge 0}-\sum_{s\ge -h'}\sum_{n< 0}\Bigr)
(-1)^n\,\tr_{\F_{m+kn,-l-2+(k+2)(n-2s)}}(\mathcal{O}_1)
\nn\\
&=\!\!&\Bigl(\sum_{n\ge 0\atop s\le h}-\sum_{n< 0\atop s> h}
 -\sum_{n\ge 0\atop s>n+h'}+\sum_{n< 0\atop s\le n+h'}\Bigr)
(-1)^n\,\tr_{\F_{m+kn,l+(k+2)(n-2s)}}(\mathcal{O}_0) 
\nn\\
&=\!\!&\Bigl(
  \sum_{-n_1-h'\leq n_2\leq n_1+h}-\sum_{n_1+h<n_2<-n_1-h'}\Bigr)(-1)^{2n_1}\,
  \tr_{\F_{m+2n_2k,l+2n_1(k+2)}}(\mathcal{O}_0),
  \nn
\ena
where the last sum is taken over $n_1,n_2\in\frac12\Z$ such that
$n_1-n_2\in\Z$. 
Note that the result is independent of the choice of $h,h'$. 

Accordingly, the sum of the zero-mode contributions coming from 
different Fock spaces can be expressed by the function 
\bea
&&
\Gamma^{(h,h')}_{m,l}(y_1,y_2|\tau) \quad(l\equiv m\bmod 2)
\lb{Gam}\\
&=\!\!&\frac{1}{\eta(\tau)^2}
  \Bigl(
  \sum_{-n_1-h'\leq n_2\leq n_1+h}-\sum_{n_1+h<n_2<-n_1-h'}\Bigr)(-1)^{2n_1}
\nn\\
&&\quad \times
e^{2\pi i \tau\bigl(
\frac{(l+1+2n_1(k+2))^2}{4(k+2)}-\frac{(m+2n_2k)^2}{4k}
\bigr)}
e^{\pi i(l+1+2n_1(k+2))y_1-\pi i (m+2n_2 k)y_2},
\nn
\ena
where 
the sum is taken over $n_1,n_2\in\frac12\Z$ with $n_1-n_2\in\Z$,
and we set $\Gamma^{(h,h')}_{m,l}(y_1,y_2|\tau)=0$ for $l\not\equiv m\bmod 2$.
Choosing $h=h'=0$, for example, we obtain the formula 
\eqref{Qlat} for the $Q_a(m)$. 
We note the properties
\bea
&&\Gamma^{(h,h')}_{m+2k,l}(y_1,y_2|\tau)
=\Gamma^{(h+1,h'-1)}_{m,l}(y_1,y_2|\tau), 
\lb{Gamsym1}\\
&&
\Gamma^{(h,h')}_{m+k,k-l}(-y_1,y_2|\tau)
=\Gamma^{(-h'+\ve',-h+\ve)}_{m,l}(y_1,y_2|\tau),
\lb{Gamsym2}
\ena
where $\ve=\bigl\{{0\;\;\;(h\not\in\Z)\atop -1\;(h\in\Z)}$ and
$\ve'=\bigl\{{1\;(h'\not\in\Z)\atop 0\;(h'\in\Z)}$.
\medskip

\noindent{\bf Modular transform}\quad
To study the continuous limit, 
it is useful to consider the modular transformation 
$\tau\rightarrow -1/\tau$. 
Unfortunately the functions $\Gamma^{(h,h')}_{m,l}(y_1,y_2|\tau)$ 
do not enjoy simple transformation properties individually.   
A way around this difficulty is proposed in \cite{Jay90}. 
Recall that up to an overall scalar factor the $Q^{(n,n)}_{l+1}(m)$ 
has the form 
\bea
\sum_{\mu,\nu\in\Z} 
f^{(n)}_{\mu,\nu}\Gamma^{(h,h')}_{m,l}(y_1,y_2|\tau)\, ,
\lb{sum}
\ena
where $f^{(n)}_{\mu,\nu}$ 
is a function of $\mu,\nu\in\Z$ given below \eqref{fn},  
and $y_1=\frac{\tau}{2k}(\mu-\nu)$, 
$y_2=\frac{\tau}{k}(\frac{2v}{k}-\frac{\mu+\nu}{2})$ 
with $v=\sum_{i=1}^n(\frac{k}{2}-v_i-u_i)$.
Making use of the independence of \eqref{sum} on $h,h'$,  
it is shown in \cite{Jay90} that 
a certain generating function of \eqref{sum} 
with respect to $m,h,h'$ can be reexpressed in terms of 
theta functions and (the zero mode contribution of) the characters of 
the $N=2$ superconformal algebra. 
{}From the knowledge of the modular property for the latter, 
the following relation can be deduced for the sum \eqref{sum}. 
\be
&&
\sum_{\mu,\nu\in\Z}
f^{(n)}_{\mu,\nu}\Gamma^{(h,h')}_{m,l}(y_1,y_2|\tau)
\\
&&
=\frac{1}{\sqrt{k(k+2)}}
\sum_{\mu,\nu\in\Z}
f^{(n)}_{\mu,\nu}
e^{-\frac{i\pi}{2\tau}\bigl((k+2)y_1^2-ky_2^2\bigr)}
\\
&&\qquad
\times \sum_{l'=0}^{k}\sum_{m'=-k}^{k-1}
\sin\frac{\pi(l+1)(l'+1)}{k+2}\,
e^{-\frac{i\pi mm'}{k}}\,
\Gamma^{(h,h')}_{m',l'}
\Bigl(\frac{y_1}{\tau},\frac{y_2}{\tau}\Bigl|-\frac{1}{\tau}\Bigr)\,.
\en
By noting \eqref{Gamsym2} and the property 
$f^{(n)}_{\nu,\mu}=f^{(n)}_{\mu,\nu}$ (see below),  
the sums in the last line can also be written as 
\be
2\sum_{l'=0}^{k}\sum_{m'=0}^{k-1} 
\sin\frac{\pi(l+1)(l'+1)}{k+2}\,
e^{-\frac{i\pi mm'}{k}}\,
\Gamma^{(h,h')}_{m',l'}
\Bigl(\frac{y_1}{\tau},\frac{y_2}{\tau}\Bigl|-\frac{1}{\tau}\Bigr).  
\en
\medskip

\noindent{\bf An identity}\quad
Let us verify the property \eqref{physical}.  
For $\mu,\nu\in\Z$, set
\bea
  &&f^{(n)}_{\mu,\nu}(u_1,\cdots,u_n;v_1,\cdots,v_n)
  \lb{fn}\\
  &=\!\!&
  \sum\,\prod_{i=1}^n\mu_i\nu_i\,\cdot\!\!
  \prod_{1\le i,j\le n}\![u_i+v_j+\sfrac{\mu_i+\nu_j}{2}]^*
  \nn\\
  &&\qquad\times
  \prod_{1\le i<j\le n}\!\!
  \frac{[u_i-u_j-\frac{\mu_i-\mu_j}{2}]^*}{[u_i-u_j]^*}
  \frac{[v_i-v_j-\frac{\nu_i-\nu_j}{2}]^*}{[v_i-v_j]^*}\,.
  \nn
\ena
Here the sum is taken over $\mu_i,\nu_i=\pm 1$ ($i=1,\cdots,n$) satisfying 
$\sum_{1=1}^n\mu_i=\mu$, $\sum_{i=1}^n\nu_i=\nu$.
It is easy to see that this function is
holomorphic and symmetric in $(u_1,\cdots,u_n)$ (resp. $(v_1,\cdots,v_n)$).
The relation \eqref{physical} reduces to the identity 
\bea
\quad\quad f^{(n)}_{\mu,\nu}(u_1,\cdots,u_n;v_1,\cdots,v_n)
=
f^{(n)}_{\nu,\mu}(u_1,\cdots,u_n;v_1,\cdots,v_n).
\lb{fnn}
\ena

Clearly \eqref{fnn} is true for $n=1$.
To show it in general,
let $g^{(n)}_{\mu,\nu}$ stand for the difference of the left hand side and
the right hand side of \eqref{fnn}.
Then
\be
  &&g^{(n)}_{\mu,\nu}(u_1,\cdots,u_n;v_1,\cdots,v_n)
  =-g^{(n)}_{\mu,\nu}(v_1,\cdots,v_n;u_1,\cdots,u_n),\\
  &&g^{(n)}_{\mu,\nu}(u_1+k,\cdots,u_n;v_1,\cdots,v_n)
  =(-1)^ng^{(n)}_{\mu,\nu}(u_1,\cdots,u_n;v_1,\cdots,v_n),\\
  &&g^{(n)}_{\mu,\nu}(u_1+\sfrac{k}{\tau},\cdots,u_n;v_1,\cdots,v_n)\\
  &&\quad=\bigl(-e^{\frac{i\pi}{\tau}}\bigr)^n
  e^{\frac{2\pi i}{k}(nu_1+v_1+\cdots+v_n+\frac{\mu+\nu}{2})}
  g^{(n)}_{\mu,\nu}(u_1,\cdots,u_n;v_1,\cdots,v_n),\\
  &&g^{(n)}_{\mu,\nu}(u_1,\cdots,u_{n-1},u;v_1,\cdots,v_{n-1},-u)\\
  &&\quad=\sum_{\mu_n=\pm 1}[\mu_n]^*
  \prod_{i=1}^{n-1}\,[u_i-u+\mu_n]^*[v_i+u+\mu_n]^*\\
  &&\qquad\qquad\times
  g^{(n-1)}_{\mu-\mu_n,\nu-\mu_n}(u_1,\cdots,u_{n-1};v_1,\cdots,v_{n-1}).
\en
{}From these properties, we conclude that $g^{(n)}_{\mu,\nu}=0$
by induction on $n$.

%%%%%%%%%%%%%%%%%%%%%%%%%%%%%%%%%%%%%%%%%%%%%%%%%%%%%%%%%%%%%%%
%                                                             %
%  References                                                 %
%                                                             %
%%%%%%%%%%%%%%%%%%%%%%%%%%%%%%%%%%%%%%%%%%%%%%%%%%%%%%%%%%%%%%%
%\bibliographystyle{unsrt}
%\bibliography{q}

\end{document}